\documentclass{amsart}
\usepackage{amsmath,amsfonts,latexsym,amssymb,amsthm}
\usepackage[mathscr]{eucal} 

\newtheorem{theorem}{Theorem}[section]
\newtheorem{corollary}[theorem]{Corollary}
\newtheorem{proposition}[theorem]{Proposition}
\newtheorem{lemma}[theorem]{Lemma}

\theoremstyle{definition}
\newtheorem{definition}[theorem]{Definition}

\newtheorem{remark}[theorem]{Remark}
\newtheorem*{claim}{Claim}
\newtheorem*{question}{Question}

\theoremstyle{definition}
\newtheorem{example}{Example}

\theoremstyle{definition}

\theoremstyle{remark}
\newtheorem*{rmk}{Remark}

\numberwithin{equation}{section}

\newcount\decno
\def\newdec{\the\decno.}
\def\newresult{\global\advance \decno by 1 {\newdec}}
\newcommand{\be}{\begin{enumerate}}
\newcommand{\ee}{\end{enumerate}}
\newcommand{\bi}{\begin{itemize}}
\newcommand{\ei}{\end{itemize}}
\newcommand{\seq}[1]{\left\langle#1\right\rangle}

\def\R{{\text{$\mathbb{R}$} }}
\def\R{{\text{$\mathbb{R}$} }}

\newcommand{\n}{{\mathfrak{n}}}
\newcommand{\m}{{\mathfrak{m}}}
\newcommand{\Lng}{\ensuremath{\mathcal{L}}}
\newcommand{\en}{\text{$\in$}}
\newcommand{\shrp}{{\text{\scriptsize\#}}}
\newcommand{\la}{\left\langle}
\newcommand{\ra}{\right\rangle}
\newcommand{\ul}{\underline}
\def\pow{{\mathcal{P}}}
\def\Lr{{\text{$L({\R})$} }}
\def\Kr{{\text{$K({\R})$} }}
\newcommand{\overcore}{\ensuremath{\overline{\core}}}
\newcommand{\overkore}{\ensuremath{\overline{\kore}}}
\newcommand{\overmouse}{\ensuremath{\overline{\mouse}}}
\newcommand{\mouse}{\ensuremath{\mathcal{M}} }
\newcommand{\nouse}{\ensuremath{\mathcal{N}} }

\def\oversystem{\text{$\langle {\langle \overline{\mouse}_\alpha \rangle}_{\alpha \in
\textup{OR}} ,  \langle \overline{\pi}_{\alpha\beta}\colon \overline{\mouse}_\alpha
\updownmap{cofinal}{1}\overline{\mouse}_\beta \rangle_{\alpha\leq\beta\in\textup{OR}}\rangle$}}
\def\prenousesystem{\text{$\langle {\langle \nouse_\gamma \rangle}_{\gamma \in
\textup{OR}} ,  \langle \pi_{\gamma}\colon \nouse
\updownmap{cofinal}{1}\nouse_\gamma \rangle_{\gamma\in\textup{OR}}\rangle$}}
\def\premousesystem{\text{$\langle {\langle \mouse_\gamma \rangle}_{\gamma \in
\textup{OR}} ,  \langle \pi_{\gamma}\colon \mouse
\updownmap{}{1}\mouse_\gamma \rangle_{\gamma\in\textup{OR}}\rangle$}}
\def\mousesystem{\text{$\langle {\langle \mouse_\alpha \rangle}_{\alpha \in
\textup{OR}} ,  \langle \pi_{\alpha\beta}\colon \mouse_\alpha \mapsigma{{n+1}}
\mouse_\beta \rangle_{\alpha\leq\beta\in\textup{OR}}\rangle$}}
\newcommand{\premouseiteration}[3]
{\text{$\langle {\langle {#1}_{#2} \rangle}_{{#2} \in \textup{OR}} , \langle
\pi_{{#2}{#3}}\colon {#1}_{#2} \updownmap{cofinal}{1}{#1}_{#3}
\rangle_{{#2}\leq{#3}\in\textup{OR}}\rangle$}}
\newcommand{\overcoresystem}{\text{$\langle {\langle \overcore_\alpha \rangle}_{\alpha \in
\textup{OR}} ,  \langle \overline{\pi}_{\alpha\beta}\colon \overcore_\alpha
\updownmap{}{1}\overcore_\beta \rangle_{\alpha\leq\beta\in\textup{OR}}\rangle$}}
\newcommand{\coresystem}{\text{$\langle {\langle \core_\alpha \rangle}_{\alpha \in
\textup{OR}} ,  \langle \pi_{\alpha\beta}\colon \core_\alpha \mapsigma{{n+1}}
\core_\beta \rangle_{\alpha\leq\beta\in\textup{OR}}\rangle$}}

\renewcommand{\iff}{\textup{ \ iff \ }}
\def\core{{\text{$\mathfrak{C}$} }}
\def\kore{{\text{$\mathcal{K}$} }}
\def\C{{\mathcal{C}}}
\def\H{{\mathcal{H}}}
\def\O{{\mathcal{O}}}
\def\AD{\textup{AD}}
\def\AC{\textup{AC}}
\def\DC{\textup{DC}}
\def\ZF{\textup{ZF}}

\def\ZFC{\textup{ZFC}}
\def\dom{{\textup{dom}}}

\def\rud{{\textup{rud}}}

\def\OR{{\textup{OR}}}
\def\Hull{{\textup{Hull}}}
\def\e{\text{$\in$}}

\newcommand{\mapsigma}[1]{\xrightarrow[\text{ \ \ $\Sigma_{#1}$}]{}}
\newcommand{\maps}[1]{\xrightarrow{\textup{#1}}}
\newcommand{\updownmap}[2]{\xrightarrow[\text{ \ \ $\Sigma_{#2}$}]{\textup{#1}}} 
\def\Jgamma{\text{$J_\gamma^{\mouse}({\R})$}}            

\def\-{\!-\!}
\def\eq{\!=\!}
\newcommand{\boldface}[2]{%
\protect\raisebox{0pt}[0pt][0pt]{%
$\underset{\displaystyle\widetilde{}}{\boldsymbol{#1}}{_{#2}}$}\mbox{\hskip 1pt}}
\newcommand{\abs}[1]{\lvert#1\rvert}
\newcommand{\skolemdefn}{3.4}
\newcommand{\deftwo}{3.88}
\newcommand{\skolem}{3.72}
\newcommand{\realpremice}{3.2}
\newcommand{\realmice}{3.4}
\newcommand{\coremice}{3.4.2}
\newcommand{\finestructure}{3}
\newcommand{\EOEL}{3.64}
\newcommand{\intismi}{3.89}
\newcommand{\skolemdef}{3.8}
\newcommand{\criterion}{3.93}
\newcommand{\soundnplusone}{3.71}
\newcommand{\newthm}{6.1}
\newcommand{\Preliminaries}{1.1}
\newcommand{\abovereals}{3.1}
\newcommand{\realonemice}{3.3}
\newcommand{\skolemimage}{3.10}
\newcommand{\acceptable}{3.15}
\newcommand{\trans}{3.18}
\newcommand{\defproj}{3.26}
\newcommand{\gammaproj}{3.30}
\newcommand{\defttho}{3.48}
\newcommand{\thmttho}{3.49}
\newcommand{\MET}{3.59}
\newcommand{\helper}{3.81}
\newcommand{\defone}{3.83}
\newcommand{\theendprop}{3.85}
\newcommand{\preddef}{3.86}
\newcommand{\predlemma}{3.87}
\newcommand{\GenDJ}{4.1}
\newcommand{\lemmatwo}{4.3}
\newcommand{\minicomp}{4.11}
\newcommand{\hulldef}{3.9}
\newcommand{\newcortwo}{6.4}
 
\begin{document}

\title[Scales and the fine structure of $K(\R)$. Part III]{Scales and the fine structure of $\boldsymbol{K(\pmb{\R})}$\\ {Part III: Scales of minimal complexity}}

\author{Daniel W. Cunningham}
\address{Mathematics Department,
State University of New York,
College at Buffalo,\\
1300 Elmwood Avenue,
Buffalo, NY 14222, USA}
\email{cunnindw@math.buffalostate.edu}

\begin{abstract} 
We obtain scales of minimal complexity in \ $\Kr$ \ using a Levy hierarchy and a fine structure theory for \ $\Kr$; \ that is,  we identify precisely those levels of the Levy hierarchy for \ $\Kr$ \ which possess the scale property. 
\end{abstract}
\maketitle
\section{Introduction}
In this paper we shall present a Levy hierarchy  for the inner model \
$\Kr$ \ and determine the minimal levels of this hierarchy which have the scale property.  As a consequence, we shall see that
in \ $\Kr$ \ there is a close connection between obtaining scales of minimal complexity and new \ $\Sigma_1$ \
truths about the reals. After we identify the levels of the Levy
hierarchy for \ $\Kr$ \ which possess the scale property, we will then be able to address the
following question, first asked in \cite{Part1}:
\begin{question}[Q]Given an iterable real premouse  $\mouse$  and  $n\ge 1$,  when does the pointclass
$\boldface{\Sigma}{n}(\mouse)$  have the scale property?
\end{question}
The boldface pointclass \ $\boldface{\Sigma}{n}(\mouse)$ \ consists of the sets of reals
definable over \ $\mouse$ \ by a \ $\Sigma_n$ \ formula allowing arbitrary constants from the domain of the
structure \ $\mouse$ \ to appear in such a definition. More generally, given \ $X\subseteq M$ \ the pointclass \
${\Sigma}_{n}(\mouse,X)$ \  consists of the sets of reals
definable over \ $\mouse$ \ by a \ $\Sigma_n$ \ formula allowing arbitrary constants from the set \ $X$ \ to appear in such a
definition. We write \ $\Sigma_n(\mouse)$ \ for the pointclass \ $\Sigma_n(\mouse,\emptyset)$. \ Throughout this paper, however,
we always allow the {\it set\/} of reals \ $\R$ \ to appear as a constant in our relevant languages (see 
\cite[subsections {\Preliminaries} \& {\abovereals}]{Part1}). Thus for \ $n\ge 1$, \ the pointclass \ $\Sigma_n(\mouse)$ \
is equal to the pointclass \ $\Sigma_n(\mouse,\{\R\})$; \ consequently, the pointclass \ $\Sigma_n(\mouse)$ \ is not
necessarily equal to the pointclass \ $\Sigma_n(\mouse,\R)$.

In \cite{Crcm} we introduced the Real Core Model \ $\Kr$ \  and showed that \ $\Kr$ \ is
an inner model containing the reals and definable scales beyond those in \ $\Lr$. \ To establish
our results in \cite{Crcm} on the existence of scales, we defined iterable real premice and extended the {\sl basic\/} fine
structural notions of Dodd-Jensen
\cite{DJ} to encompass iterable ``premice above the reals.'' Consequently, we were able to prove the
following result (see \cite[Theorem 4.4]{Crcm}):

\begin{theorem}\label{firstscales}  Suppose that \ $\mouse$ \  is an
iterable real premouse and that \ $\mouse \models \AD$. \ Then \ $\Sigma_1(\mouse)$ \
has the scale property. 
\end{theorem}
By allowing for real parameters in the proof of Theorem~\ref{firstscales}, we have the following corollary: 

\begin{corollary}\label{relscales}  Suppose that \ $\mouse$ \  is an
iterable real premouse and that \ $\mouse \models \AD$. \ Then \ $\Sigma_1(\mouse,\R)$ \
has the scale property. 
\end{corollary}

We say that \ $\mathcal{M}=(M,\R,\kappa,\mu)$ \ is a real 1--mouse  if \ $\mouse$ \ is an iterable real premouse and \
$\mathcal{P}(\R\times\kappa)\cap \boldface{\Sigma}{1}(M)\not\subseteq M$, \ where \ $M$ \ has the form \ $J_\alpha[\mu](\R)$ \ and \
$\kappa$ \ is the ``measurable cardinal'' in \ $\mouse$ \ (see \cite[subsections {\realpremice} \& {\realonemice}]{Part1}). Real
1--mice suffice to define the real core model and to prove the results in \cite{Crcm} about \ $\Kr$; \ however, real 1--mice are not
sufficient to construct scales of minimal complexity in \ $\Kr$. \ Our solution to the problem of identifying these scales
requires the development of a {\sl full\/} fine structure theory for \ $\Kr$. \ In the paper \cite{Cfsrm} we initiated this
development by generalizing Dodd-Jensen's notion of a mouse to that of a {\it real mouse\/} (see \cite[subsection
{\realmice}]{Part1}). This is accomplished by replacing \ $\Sigma_1$ \ with \ $\Sigma_n$, \ where \ $n$ \ is the smallest
integer such that \
$\mathcal{P}(\R\times\kappa)\cap\boldface{\Sigma}{n+1}(\mouse)\not\subseteq M$, \ together with a stronger iterability condition.

Let \ $\mouse$ \ be a real mouse and assume that there is an integer \ $m\ge 1$ \
such that \ $\mathcal{P}(\R)\cap\boldface{\Sigma}{m}(\mouse)\not\subseteq M$. \ We shall let
\ $\m=m(\mouse)$ \ denote the least such integer and, in this case, we say that \ $\mouse$ \  is \
{\it weak} \ if 
\be
\item[(1)] $\mouse$ \ is a proper initial segment of an iterable real premouse, and 
\item[(2)] $\mouse$ \ realizes a \ $\Sigma_\m$ \ type not realized in any proper initial segment of \ $\mouse$. 
\ee
In (2) a \ $\Sigma_\m$ \ type is a non-empty subset \ $\Sigma$ \ of \[\{ \theta\in \Sigma_\m\cup\Pi_\m  : \text{$\theta$
is a formula of one free variable}\}\] and \ $\mouse$ \ is said to realize \ $\Sigma$ \ if there is an \ $a\in M$ \
such that \ $\mouse\models\theta(a)$ \ for all \ $\theta\in\Sigma$. \ Finally, if \ $\mouse$ \ satisfies (1) but fails to
satisfy (2), then we say that \ $\mouse$ \ is {\it strong\/}.

Using the fine structure of real mice developed in \cite{Cfsrm} and
\cite{Part1}, we establish in Part II \cite[Theorem~{\newthm}]{Part2}
the following theorem on the existence of scales:
\begin{theorem}\label{newthmrpt}  Suppose that  \ $\mouse$ \  is a
weak real mouse satisfying \ $\AD$. \ Then \
$\boldface{\Sigma}{\m}(\mouse)$ \ has the scale property, where \ $\m=m(\mouse)$. 
\end{theorem}

\begin{remark} We specifically note that the proofs of both Theorem~\ref{firstscales} and Theorem~\ref{newthmrpt}  require {\it
only\/} the determinacy of sets of reals in the relevant \ $\mouse$. \ This is a critical property for proving that \
$\Kr$ \ satisfies \ $\AD$ \ under certain hypotheses.  
\end{remark}

Theorem~\ref{newthmrpt} is an essential component in our analysis of scales in \ $\Kr$. \ Our next result, which
follows from the proof of Theorem~\ref{newthmrpt}  (see \cite[Theorem~{\newcortwo}]{Part2}), is also an important ingredient in
our analysis of scales in \ $\Kr$.
\begin{theorem}\label{essential} Suppose that  \ $\mouse$ \  is a
weak real mouse and let \ $\m=m(\mouse)$. \ For any set of reals \ 
$P\in \boldface{\Sigma}{\m}(\mouse)$, \ there exists a total \ $\boldface{\Sigma}{\m}(\mouse)$ \ 
map \ $k\colon\omega\to M$ \ such that \ 
$P=\bigcup_{i\in\omega}k(i)$. 
\end{theorem}

Theorems~\ref{firstscales} and \ref{newthmrpt} are the key results that we will use in
this paper to give a complete description of those levels of the Levy hierarchy for \ $\Kr$ \ which
have the scale property.  The work presented here can be viewed as a generalization of Steel's work on the existence of
scales in the inner model \ $\Lr$. \ Steel \cite{Steel} develops a fine structure theory and a Levy
hierarchy for \ $\Lr$. \ Using this development, Steel solves the problem of finding scales of
minimal complexity in \ $\Lr$ \ and, as a consequence, shows that there is a close connection between
obtaining such scales and new \ $\Sigma_1$ \ truths in \ $\Lr$ \ about the reals.

Our paper is organized into \ref{questionQ} sections. In Section~\ref{realcm} we review the definitions of \
$\Lr$ \ and \ $\Kr$. \ Section~\ref{gaps} discusses the notion of a $\Sigma_1$--gap. Section~\ref{minimal} presents a complete
description of those levels of the Levy hierarchy for iterable real premice which have the scale property. 
Section~\ref{extensions} shows that the premouse iteration of a premouse preserves its $\Sigma_1$--gaps and preserves its
internal pointclasses (see Definition~\ref{internal}). We define a Levy hierarchy for \ $\Kr$ \ in Section
\ref{minimalKr} and then construct scales of minimal complexity in \ $\Kr$. \
In Section~\ref{questionQ} we present an answer to the question posed at the beginning of this paper.  

\subsection{Preliminaries and notation}\label{Prelims}

Let \ $\omega$ \ be the set of all natural numbers.  \ $\R = {^\omega \omega}$ \ is
the set of all functions from \ $\omega$ \ to \ $\omega$. \ We call \ $\R$ \ the set of
reals and regard \ $\R$ \ as a topological space by giving it the product topology,
using the discrete topology on \ $\omega$. \ For a set \  $A \subseteq \R$ \ we
associate a two person infinite game on \ $\omega$, \ with {\it payoff} \
$A$,
\ denoted by \ $G_A$:
\[
\begin{aligned}[c]
{}&{\mathbf{I}}\phantom{{\mathbf{I}}} \qquad x(0) \qquad \phantom{x(1)} \qquad x(2) \qquad \phantom{x(3)}\quad \\
{}&{\mathbf{I}}\mathbf{I} \qquad \phantom{x(0)}\qquad x(1) \qquad \phantom{x(2)} \qquad x(3) \quad 
\end{aligned}
\begin{gathered}[c]
{\cdots}
\end{gathered}
\] 
in which player \ {\bf I} \ wins if \ $x \in A$, \ and \ {\bf II} \ wins if \ $x
\notin A$. \ We say that \ $A$ \ is {\it determined\/} if the corresponding game \ $G_A$ \
is determined, that is, either player {\bf I} or  {\bf II} has a winning strategy (see \cite[page
287]{Mosch}).  The {\it axiom of determinacy\/} \ ($\AD$) \ is a regularity hypothesis about games on \
$\omega$ \ and states: \ $\forall A \subseteq \R \  ( A  \text{\ is determined})$.

We work in \  $\ZF+\DC$ \ and state any additional hypotheses as we need them, to keep a close
watch on the use of determinacy in the proofs of our main theorems.  Variables \ $x, y, z, w \dots$
\ generally range over \ $\R$, \ while \ $\alpha, \beta, \gamma, \delta \dots$ \ range over \ $\OR$, \ the
class of ordinals.  The cardinal \ $\Theta$ \ is the supremum of the ordinals which are the surjective image
of \ $\R$. 

A {\it pointclass\/} is a set of subsets of \ $\R$ \ closed under recursive substitutions. A
boldface pointclass is a pointclass closed under continuous substitutions.  For a
pointclass \ $\Gamma$, \ one usually writes \ ``$\Gamma\-\AD$'' or ``Det($\Gamma$)'' \ to denote
the assertion that all games on \ $\omega$ \ with payoff in \ $\Gamma$ \ are
determined. For the notions of a scale and of the scale property (and any other
notions from Descriptive Set Theory which we have not defined), we refer the reader
to Moschovakis \cite{Mosch}.

A proper class \ $M$ \ is called an {\it inner model\/} if and only if \ $M$ \ is a
transitive  $\in$--model of \ $\ZF$ \ containing all the ordinals. For an inner model \ $M$ \ with \
$X\in M$ \ we shall write \ $\mathcal{P}^{M}(X)$ \ to denote the power set of \ $X$ \ as computed in
\ $M$. \  For an ordinal \ $\kappa\in M$, \ we shall abuse
standard notation slightly and write \ ${^\kappa M}=\{f\in M \ \vert \ f\colon \kappa \rightarrow M \}$. \ 

We distinguish between the notations \ $L[A]$ and $L(A)$. \ 
The inner model \ $L(A)$ \ is defined to be the class of sets constructible {\it above\/} \ $A$, \ that is, \ one starts
with a set \ $A$ \ and iterates definability in the language of set theory. \ Thus, \ $L(A)$ \ is  the smallest
inner model \ $M$ \ such that \ $A\in M$. \ The inner model \ $L[A]$ \ is defined to be the class of sets constructible {\it
relative\/} to \ $A$, \ that is, \ one starts with the empty set and iterates definability in the language of set theory
augmented by the predicate \ $A$. \ Consequently, \ $L[A]$ \ is the smallest inner model \ $M$ \ such that \ $A\cap M\in
M$ \ (see page 34 of \cite{Kana}). Furthermore, one defines \ $L[A,B]$ \ to be the class of sets constructible {\it relative\/}
to \ $A$ \ and \ $B$, \ whereas \ $L[A](B)$ \ is defined as the class of sets constructible {\it relative\/} to \ $A$ \ and 
{\it above\/}
\ $B$. \ Thus, \ $A\cap L[A](B)\in L[A](B)$ \ and \ $B\in L[A](B)$.

Given a model \ $\mouse = (M,c_1,c_2,\dots,c_m,A_1,A_2,\dots,A_N),$ \ where the \  $A_i$ \
are predicates and the \ $c_i$ \ are constants,  if \ $X\subseteq M$ \ then \
$\Sigma_n(\mouse,X)$ \ is the class of relations on \ $M$ \ definable over \ $\mouse$ \ by a
\ $\Sigma_n$ \ formula from parameters in \ $X \cup \lbrace c_1,c_2,\dots,c_m\rbrace$. \
$\Sigma_\omega(\mouse,X) =  \bigcup_{n\in\omega} \Sigma_n(\mouse,X).$  We write
``$\Sigma_n(\mouse)$'' for \ $\Sigma_n(\mouse,\emptyset)$ \ and \
``$\boldface{\Sigma}{n}(\mouse)$'' for  the boldface class \ $\Sigma_n(\mouse,M).$ \ 
Similar conventions hold for \ $\Pi_n$ \ and \ $\Delta_n$ \ notations.  If \ $\mouse$ \ is
a substructure of \ $\nouse$ \ and \ $X\subseteq M \subseteq N$, \ then 
``$\mouse \prec_n^X \nouse$'' means that \ $\mouse \models \phi[a]$ \ if and only if \ $\nouse \models \phi[a]$, \ for all \ $a
\in (X)^{<\omega}$ \ and for all \
$\Sigma_n$ \ formulae \ $\phi$ \ (the formula \ $\phi$ \ is allowed constants taken
from \ $\lbrace c_1,c_2,\dots,c_m\rbrace$). \ We write \ ``$\mouse \prec_n \nouse$'' \ 
for \ ``$\mouse \prec_n^M \nouse$.'' \ In addition, for any two models \ $\mouse$ \ and \
$\nouse$, \  we write \ $\pi:\mouse\mapsigma{n}\nouse$ \ to indicate that the map \ $\pi$ \ is
a \ $\Sigma_n$--elementary embedding, that is, \ $\mouse \models \phi[a]$ \ if and only if \ $\nouse \models \phi[\pi(a)]$, \
for all \ $a=\langle a_0,a_1,\dots\rangle \in (M)^{<\omega}$ \ and for all \ $\Sigma_n$ \ formulae \ $\phi$, \ where \ $0\le
n\le\omega$ \ and \ 
$\pi(a)=\langle \pi(a_0),\pi(a_1),\dots\rangle$.

We now give a brief definition of an iterable real premice (for more details
see \cite{Part1} or \cite{Crcm}). First, we present a preliminary definition.

\begin{definition}\label{def:Rcomplete} Let \ $\mu$ \ be a normal measure on \
$\kappa$. \ We say that \ $\mu$ \ is an {\it $\R$--complete measure} on \ $\kappa$ \ 
if the following holds: if \ $\langle A_x : x\in \R\rangle$ \ is any sequence such that \ 
$A_x\in \mu$ \ for all \ $x\in \R$, \ then \ $\bigcap_{x\in \R}A_x \in \mu$.
\end{definition}

For \ $N\in\omega$, \ the language 
\[\Lng_N=\lbrace\en,\underline{\R},\underline{\kappa},\mu,A_1,\dots,A_N\rbrace\]
consists of the constant symbols \ $\underline{\R}$ \ and \
$\underline{\kappa}$ \ together with the membership relation \ $\in$ \ and the predicate symbols \ $\mu, A_1,\dots,A_N$.

\begin{definition}\label{premice} A model \ $\mouse = 
(M,\en,\underline{\R}^{\mouse},\underline{\kappa}^{\mouse},\mu,A_1,\dots,A_N)$ \ is
a \ {\it premouse} (above the reals) \ if 
\be
\item $\mouse$ \ is a transitive set model of \ $V=L[\mu,A_1,\dots,A_N](\underline{\R})$
\item $\mouse \models \text{``$\mu$ is an \ $\underline{\R}$--complete \ measure on 
\ $\underline{\kappa}$''}$. 
\ee
$\mouse$ \ is a \ {\it pure premouse} \ if \ $\mouse = 
(M,\underline{\R}^{\mouse},\underline{\kappa}^{\mouse},\mu)$ \ (that is, $N=0$). \ Finally,  \ $\mouse$ \  is a \ {\it
real premouse} \  if it is pure and \ ${\underline{\R}}^{\mouse}= \R$.
\end{definition}

A real premouse \ $\mouse$ \ has a natural Jensen hierarchy. For any \
$\alpha \in \OR^{\mouse}$ \ we let \  $S_\alpha^{\mouse}({\R})$ \ denote the unique set
in \ $\mouse$ \  satisfying  \ $\mathcal{M} \models \exists f(\varphi(f) \land \alpha\in \dom(f)
\land  S_\alpha^{\mouse}({\R}) = f(\alpha))$, \  
where \ $\varphi$ \ is the \ $\Sigma_0$ \ sentence used to define the sequence \ 
$\langle S_\gamma^{\mouse}({\R})\ : \ \gamma <
\OR^{\mouse} \rangle$ \ (see Definition 1.5 of \cite{Crcm}).
For \ $\lambda = \OR^{\mouse}$, \ let \ 
$S_\lambda^{\mouse}({\R})  = \bigcup\limits_{\alpha < \lambda} S_\alpha^{\mouse}({\R})$.
\   Let \  $\widehat{\OR}$ \ denote the class of ordinals \ $\lbrace \gamma : \text{ the ordinal }
\omega \gamma \text{ exists} \rbrace$ \ and let \ 
$\Jgamma = S_{\omega\gamma}^{\mouse}({\R})$, \ for \ $\gamma \leq
\widehat{\OR}^{\mouse}$. \ It follows that
\  $M=J_\alpha^{\mouse}({\R})$ \ where \ $\alpha=\widehat{\OR}^{\mouse}$. \ 
Let \ $\mouse^\gamma$ \ be the
substructure of \ $\mouse$ \ defined by \ $\mouse^\gamma =
(\Jgamma,\en,{\R},\Jgamma \cap \mu)$ \  
for \ $1<\gamma \leq \widehat{\OR}^{\mouse},$ \ and let \ 
$M^\gamma = \Jgamma$. \ We can write \
${\mouse}^\gamma = (M^\gamma,\en,{\R},\mu),$ \ as this will
cause no confusion. In particular, \ $\mouse^\gamma$ \ is {\it amenable}, that is, \ 
$a\cap \mu \in M^\gamma$, \ for all \ $a\in M^\gamma$.

Given a premouse \ $\mouse$ \ we can construct iterated ultrapowers and obtain a
commutative system of models by taking direct limits at limit ordinals.
\begin{definition}\label{pmiteration}  Let  \ $\mouse$ \  be a premouse. Then 
\begin{equation}\premousesystem\label{lablepmiteration}\end{equation}
is the commutative system satisfying the inductive definition:
\be
\item $\mouse_0 = \mouse$
\item $\pi_{\gamma\gamma} = \text{identity map, \ and } \pi_{\beta\gamma} \circ
\pi_{\alpha\beta} = \pi_{\alpha\gamma} \text{ \ for all } 
\alpha\leq\beta\leq\gamma\leq\lambda$
\item If \ $\lambda=\lambda^\prime + 1$, \ then \ $\mouse_\lambda =
\text{ultrapower of } \mouse_{\lambda^\prime}$, \ and \ 
$\pi_{\alpha\lambda}=\pi^{\mouse_{\lambda^\prime}} \circ
\pi_{\alpha{\lambda^\prime}}$ \ for all \  $\alpha\leq\lambda^\prime$
\item If \ $\lambda$ \ is a limit ordinal, \ then \ $\langle\mouse_\lambda,
\langle \pi_{\alpha\lambda}\colon \mouse_\alpha \rightarrow \mouse_\lambda
\rangle_{\alpha<\lambda}\rangle$ \ is the direct limit of \ 
\[\la {\la \mouse_\alpha \ra}_{\alpha < \lambda} , 
\la \pi_{\alpha\beta}\colon \mouse_\alpha \rightarrow \mouse_\beta
\ra_{\alpha\leq\beta<\lambda}\ra.\]
\ee
\end{definition}
The commutative system in the above (\ref{lablepmiteration}) is called the  {\it premouse iteration} of \
$\mouse$. \  We note that the maps in the above commutative system are cofinal and 
are \ $\Sigma_1$ \ embeddings, that is,
\[\pi_{\alpha\beta}\colon \mouse_\alpha \updownmap{cofinal}{1} \mouse_\beta\]
for all \ $\alpha\leq\beta\in \OR$. \ We shall call \  $\pi_{0\beta}\colon \mouse\mapsigma{1}
\mouse_\beta$ \ the {\it premouse embedding\/} of \ $\mouse$ \ into its \ $\beta^{\,\ul{\text{th}}}$ \
{\it premouse iterate\/} \ $\mouse_\beta$.
\begin{definition} A premouse  \ $\mouse$ \  is an {\it iterable
premouse} if  \ $\mouse_\lambda$ \ is well-founded for all \ $\lambda \in \OR$.
\end{definition}

For an iterable premouse  \ $\mouse$ \  and \ $\alpha\in \OR$, \ we identify \
$\mouse_\alpha$ \ with its transitive collapse. It follows that
\[\mouse_\alpha = (M_\alpha,\en,{\underline\R}^{\mouse_\alpha},
\underline\kappa^{\mouse_\alpha},\mu,A_1,\dots,A_N)\] 
is a premouse and we write
\ $\kappa_\alpha = \pi_{0\alpha}(\underline{\kappa}^{\mouse}) = 
\underline{\kappa}^{\mouse_\alpha}$ \ for \ $\alpha\in \OR$.

Whenever we write \ $\boldface{\Sigma}{n}(\mouse)=\boldface{\Sigma}{n}(\nouse)$, \ we implicitly mean that this is an
equality between pointclasses. 

\begin{definition}Let \ $\mouse$ \ and \ $\nouse$ \ be real premice. We shall say that \
$\boldface{\Sigma}{n}(\mouse)=\boldface{\Sigma}{n}(\nouse)$ \ {\it as pointclasses\/} if for every \ $\Sigma_n$ \ formula of one
variable \ $\varphi(v)$ \ with constants from \
$\mouse$, \ there is a \ $\Sigma_n$ \ formula of one variable \ $\psi(v)$ \ with constants from \ $\nouse$ \ such
that \ $\mouse\models\varphi(x)$ \ if and only if \ $\nouse\models\psi(x)$ \ for all \ $x\in\R$, and vice versa. 
\end{definition}

We now recall the definition of \ $\Sigma_\omega^\mu$ \ formulae and some other notions from \cite{Part1}. A real premouse \
$\nouse=(N,\R,\kappa,\mu)$  \ is a model of the language \
$\Lng=\{\in,\underline{\R},\underline{\kappa},\mu\}$ \ where \ $\mu$ \ is a predicate. At times, we will want to add a quantifier
to the language \ $\Lng$. \ Since the quantifier extends the predicate \ $\mu$ \ in our intended structures, we shall use the same
symbol \ $\mu$ \ for this quantifier. We shall denote this expanded language by \ $\Lng^\mu$ \ and write \
$\Sigma_\omega^\mu$ \ for the formulae in this expanded language. For \ $\gamma$ \ such that \
$\kappa<\gamma<\widehat{\OR}^{\nouse}$, \ let \ $\mu^{\gamma+1}=N^{\gamma+1}\cap\mu$. \ Then \
$(\nouse^\gamma,\mu^{\gamma+1})$ \ is an \ $\Lng^\mu$ \ structure, where the new
quantifier symbol is to be interpreted by \ $\mu^{\gamma+1}$. \ That is, \ $(\nouse^\gamma,\mu^{\gamma+1})\models
(\mu\,\varkappa\in\kappa)\psi(\alpha)$ \ if and only if \ $\{\varkappa\in\kappa : (\nouse^\gamma,\mu^{\gamma+1})\models
\psi(\varkappa)\}\in\mu^{\gamma+1}$. \ The following is Definition~{\preddef} of \cite{Part1}, but with an
additional clause.

\begin{definition}\label{preddefagain} Let \ $\nouse=(N,\R,\kappa,\mu)$ \ be a real premouse and 
let \ $\kappa<\gamma<\widehat{\OR}^{\nouse}$. \ We shall say that \ $\mu^{\gamma+1}$ \ is {\it predictable\/} if 
the following holds: For each \ $\Lng$ \ formula \ $\chi(v_0,v_1,\dots,v_k)$ \ there is another  \ $\Lng$ \
formula \ $\psi(v_1,\dots,v_k)$ \ with a parameter  \ $d\in N^\gamma$ \ such that for all \ $a_1,\dots,a_k\in N^\gamma$ \ 
\begin{equation}B_{a_1,\dots,a_k}\in\mu^{\gamma+1}\iff \nouse^\gamma\models\psi(a_1,\dots,a_k)\label{dictable}\end{equation}
where \ $B_{a_1,\dots,a_k}=\{\varkappa\in\kappa : \nouse^\gamma\models \chi(\varkappa,a_1,\dots,a_k)\}$. \ If there is single
parameter  \ $d\in N^\gamma$ \ that satisfies (\ref{dictable}) for all such formula \ $\psi$, \ then we will state that \
$\mu^{\gamma+1}$ \ is  $d$--{\it predictable\/}. 
\end{definition}

We will now go over the definition of a real 1--mouse.
\begin{definition} Let \ $\mouse$ \ be an iterable real premouse. 
The {\it projectum} \ $\rho^{}_\mouse$ \ is the least ordinal \ $\rho \le
\widehat{\OR}^{\mouse}$ \ such that \  $\pow(\R\times
\omega\rho)\cap\boldface{\Sigma}{1}(\mouse)\not\subseteq M$, \ and \ $p^{}_\mouse$ \ is the \ 
$\le_{BK}$--least \ $p\in[\OR^{\mouse}]^{<\omega}$ \ such that \  $\pow(\R\times
\omega\rho^{}_\mouse)\cap\Sigma_1(\mouse,\{p\})\not\subseteq M$. 
\end{definition}

\begin{definition}\label{onemouse}  An iterable real premouse \ $\mouse$ \ is a \ {\it real 1--mouse} \ if \
$\omega\rho^{}_{\mouse}\le \kappa^{\mouse}$. 
\end{definition}

Real 1--mice suffice to define the class \ $\Kr$ \ and to prove the results in
\cite{Crcm};  however, real 1--mice are not sufficient to construct scales of minimal complexity. Our
solution to the problem of identifying these scales in \ $\Kr$ \ requires the development of a {\sl full\/}
fine structure theory for \ $\Kr$. \ In the paper \cite{Cfsrm} we initiated this development by generalizing
Dodd-Jensen's notion of a mouse to that of a {\it real mouse\/}  (see subsection~{\realmice} of \cite{Part1}). This is
accomplished by 
\bi
\item isolating the concept of {\it acceptability above the reals\/}\footnote{This concept extends the Dodd-Jensen notion of
acceptability to include the set of reals.}  (see \cite[Definition
{\acceptable}]{Part1}),
\item replacing \ $\Sigma_1$ \ with \ $\Sigma_n$, \ where \ $n$ \ is the smallest integer such that \ 
$\mathcal{P}(\R\times\kappa)\cap\boldface{\Sigma}{n+1}(\mouse)\not\subseteq M$, 
\item defining an iteration procedure stronger than premouse iteration.
\ei

Let \ $\mouse=(M,\R,\underline{\kappa}^{\mouse},\mu)$ \  be an iterable real premouse. The \
$\Sigma_1$--{\it master code\/} \ $A_\mouse$ \ of \ $\mouse$ \  is the set 
\[A_\mouse = \{(x,s)\in \R\times (\omega\rho^{}_\mouse)^{<\omega} : \mouse\models
\varphi_{x(0)}[\lambda n.x(n+1),s,p^{}_\mouse]\}\]
where \ $\langle\varphi_i : i\in\omega\rangle$ \ is a fixed recursive listing of all the \
$\Sigma_1$ \ formulae of three variables in the language \ $\Lng=\{\,\en,\underline{\R},\underline{\kappa},\mu\,\}$.

Theorem~{\GenDJ} of \cite{Part1} proves that an iterable real premouse \ $\mouse$ \ is acceptable above the reals. Using 
\ $A_\mouse$ \ and \ $\omega\rho^{}_\mouse$ \ one defines a new structure with domain \ $H^{\mouse}_{\omega\rho^{}_\mouse}=\{a\in
M :  \abs{T_c(a)}_\mouse<\omega\rho^{}_\mouse\}$ \ where \ $T_c(a)$ \ denotes ``the transitive closure of \ $a$'' and \
$\abs{a}_M$ \ denotes the least ordinal \ $\lambda$ \ in \ $\mouse$ \ such that \ $f\colon \lambda \times
{\R}\maps{onto} a$ \ for some \ $f\in M$. \ Let \ $M^1=H^{\mouse}_{\omega\rho^{}_\mouse}$. \ The \ $\Sigma_1$-code of \ $\mouse$
\ is the structure \ $\mouse^1 = (M^1,\R,\underline{\kappa}^{\mouse},\mu, A_1)$, \ where \ $A_1=A_\mouse$. \ Because \ 
$\mouse$ \ is acceptable,  we can repeat this construction and inductively  define structures \ $\mouse^n$ \ in the language \
$\Lng_n=\lbrace\e,\underline{\R},\underline{\kappa},\mu, A_1,\dots,A_n\rbrace$ \  where the predicate symbols \
$A_1,\dots,A_n$  \ represent the previously defined master codes and thus, one can define \ $\rho^n_\mouse$ \ and \
$A_{\mouse^n}$. \  When there is an integer \ $\n$ \ such that \ $\rho_{\mouse}^{\n+1}\le \kappa^\mouse < \rho_{\mouse}^\n$, \
then we say that \ $\mouse$ \ is {\it critical\/} \ and we let
\ $\n=n(\nouse)$ \ denote this integer. If the structure \ $\mouse^\n$ \ is sufficiently
iterable, then we say that \
$\mouse$ \ is a real mouse. More specifically,   let \ $\overline{\mouse}=\mouse^\n$. \ Since \ $\overline{\mouse}$ \  is an
iterable real premouse, let
\begin{equation}\oversystem\label{prmc}\end{equation}
be the  {\it premouse iteration\/} of \ $\overline{\mouse}$ \ as in Definition~\ref{pmiteration}.
We can extend the system (\ref{prmc}) of transitive models via the extension of embeddings lemma (Lemma~{\EOEL} of
\cite{Part1}) and obtain the commutative system of transitive structures
\begin{equation}\mousesystem.\label{mouse.iter}\end{equation}
The system (\ref{mouse.iter}) is called the {\it mouse iteration\/} of \ $\mouse$. \ We shall call  \
$\pi_{0\beta}\colon \mouse\mapsigma{\n+1} \mouse_\beta$ \ the {\it mouse embedding\/} of \ $\mouse$ \ into its \
$\beta^{\,\ul{\text{th}}}$ \ {\it mouse iterate\/} \ $\mouse_\beta$.

\begin{rmk} A real 1--mouse \ $\mouse=(M,\R,\kappa,\mu)$ \ is the simplest of real mice; because \ $\mouse$ \ is
iterable and \ $\mathcal{P}(\R\times\kappa)\cap\boldface{\Sigma}{1}(M)\not\subseteq M$.
\end{rmk}

We now review the definition of the core of a real mouse
\ $\mouse$. \ Let \ $\overline{\mouse}=\mouse^\n$, \ where \ $\n=n(\mouse)$, \ and let \ $\mathcal{H} =
\Hull_1^{\overline{\mouse}}({\R\cup \omega\rho^{}_{\,\overline{\mouse}}}
\cup \{p^{}_{\,\overline{\mouse}}\})$.
\ Thus, \ $\H\prec_1\overline{\mouse}$. \ Let \ $\overcore$ \ be the transitive collapse of \ $\H$. \ By Lemma~{\EOEL} of
\cite{Part1} there is a decoding \ $\core$ \ of \ $\overcore$ \ and a map \ $\sigma\colon \core \mapsigma{{n+1}} \mouse$. \ 
It follows that \ $\core$ \ is a real mouse with \ $n(\core)=n(\mouse)$. \ We denote \ $\core$ \ by \ $\core(\mouse)$. \ Let
\[\overcoresystem\]
be the premouse iteration of \ $\overcore$ \ and let
\[\coresystem\]
be the mouse iteration of \ $\core$. It follows that \
$\mouse$ \ is a mouse iterate of \ $\core$; \ that is, there is an ordinal \ $\theta$ \ such that \ $\overcore_\theta=\overmouse$ \ and \
$\core_\theta=\mouse$.

Let \ $\kappa_\alpha = \pi_{0\alpha}(\underline{\kappa}^{\core})$ \ for \ $\alpha\in \OR$. \ In particular,
\ $\kappa_0=\underline{\kappa}^\core$. \ Given that \ $\core_\theta=\mouse$, \
let \ $I_m=\{\,\kappa_\alpha : \alpha \text{ is $m$--good} \land \alpha < \theta\,\}$, \ where an ordinal \ $\alpha$ \ is  \
$m$--good \ if and only if \ $\alpha$ \ is a multiple of \ $\omega^m$. \ Corollary~{\helper} of \cite{Part1} asserts that \
$I_m$ \ is a set of order \ $\Sigma_{m+1}(\mouse,\{\,\pi_{0\theta}(a) : a\in C\,\})$ \ indiscernibles, \ where \ $C$  \ is the
domain of \
$\core$. \ We  now state a result from \cite{Part1} that we use in subsection~\ref{endgap} when we deal with the
existence of scales at the ``end of a gap.''

\begin{lemma}\label{deftwotwo}
Let  \ $\mouse$ \ be a mouse with core \ $\core$ \ and let \ $\n=n(\mouse)$. \ Suppose that \ $\rho_{\mouse}^{\n+1}
< \kappa^{\mouse}$ \ and that \ $\mouse=\nouse^\gamma$ \ for an iterable premouse \ $\nouse$ \ where \
$\kappa^{\mouse}=\kappa^{\nouse}=\kappa$ \ $\kappa<\gamma<\widehat{\OR}^{\nouse}$. \ Let \ $\theta$
\ be such that the mouse iterate \ $\core_\theta=\mouse$ \ and let \ $\kappa_0=\kappa^\core$. \ Then 
\be
\item[(1)] $\theta=\kappa$, \ $\theta$ \ is a multiple of \ $\omega^\omega$ \ and \ $I_m\in\mu^\nouse$ \ for all \ $m\in\omega$
\item[(2)] $I_m$ \ is uniformly \ $\Sigma_\omega(\overline{\mouse},\{\kappa_0\})$
\item[(3)] $I_m$ \ is uniformly \ $\Sigma_\omega({\mouse},\{\kappa_0\})$ and its definition depends only on \
$\n$
\item[(4)] $\mu^{\gamma+1}$ \ is $\kappa_0$--predictable.
\ee
\end{lemma}
\begin{proof} Items (1) and (2) of the above list follow directly from Lemma~{\defone} and Lemma
{\deftwo}  in \cite{Part1}.  Since
each \ $I_m$ \ is uniformly \ $\Sigma_\omega(\overline{\mouse},\{\kappa_0\})$ \ and because \ $\overmouse$ \
is definable over \ $\mouse$ \ with a definition depending only on \ $\n$, \ we see that (3) holds.
Our proof of Lemma~{\deftwo} in \cite{Part1} assumes that \ $I_m$ \ is \
$\boldface{\Sigma}{\omega}(\mouse)$ \ and proves, as a claim, that \
$\mu^{\gamma+1}$ \ is predictable. However, since each \ $I_m$ \ is \
$\Sigma_\omega(\mouse,\{\kappa_0\})$, \ our proof of Lemma~{\deftwo} is easily modified to show that
\ $\mu^{\gamma+1}$ \ is $\kappa_0$--predictable.
\end{proof}

Recall that a real premouse \ $\mouse$ \ is a pure premouse; that is, it has the form \ $\mouse =
(M,\R,\underline{\kappa}^{\mouse},\mu)$. \  For the remainder of this paper, when we say that a structure \ $\nouse$ \ is a
premouse we shall mean, for the most part,  that \ $\nouse$ \ is real premouse. It will be clear from the context when we are
actually working with premice that are not pure. For more details on the matters discussed in this section, see \cite[Section
{\finestructure}]{Part1}.

\section{The Real Core Model}\label{realcm}
In this section we shall review the basic definitions of  \ $\Lr$ \ and \ $\Kr$.

\subsection{The inner model {\mathversion{bold}$L(\mathbb{R})$}}\label{jensenL(R)}
We review the Jensen hierarchy for \ $\Lr$. \ We presume the reader is familiar with the rudimentary functions (see
\cite{Jensen}). Let \ $\rud(M)$ \ be the closure of \ $M\cup\{M\}$ \ under the rudimentary functions. Let 
\begin{align*}
J_1(\R) &=V_{\omega+1} = \{a : \textup{rank}(a)\le\omega\},\\
J_{\alpha+1}(\R) &=\rud(J_{\alpha}(\R)) \ \textup{ for $\alpha>0$,}\\
J_{\lambda}(\R) &=\bigcup_{\alpha<\lambda}J_{\alpha}(\R)) \ \textup{ for limit ordinals $\lambda$}.
\end{align*}
$\Lr=\bigcup\nolimits_{\alpha\in\OR}J_{\alpha}(\R))$ \ is the smallest inner model of $\ZF$ containing the reals.
Under the hypothesis that \ $\Lr$ \ is a model of \
$\AD + \DC$, \ researchers have essentially settled all the important problems concerning the descriptive set theory and structure
of \ $\Lr$.

\subsection{The inner model {\mathversion{bold}$K(\mathbb{R})$}}
Using real 1--mice there a natural way to define an inner model larger than \ $\Lr$.
\begin{definition} The real core model is the class 
 \ $\Kr = \{\,x : \exists\, \nouse (\nouse \text{ \ is a real 1--mouse} \land x\in N)\,\}$.
\end{definition}

One can prove that \ $\Kr$ \ is an inner model of \ $\ZF$ \ and
contains a set of reals not in \ $\Lr$ \ (see \cite{Crcm}). It turns out that problems concerning the descriptive set theory and
structure of \ $\Kr$ \ can also be settled under the hypothesis that \ $\Kr$ \ is a model of \ $\AD$. \ For example, using a
mixture of descriptive set theory, fine structure and the theory of iterated ultrapowers,  one can produce definable scales in \
$\Kr$ \ beyond those in \ $\Lr$ \ and prove that \ $\Kr\models\DC$.

\begin{remark} There exists an iterable real premouse if and only if \ $\R^\shrp$ \ exists.
So \ $\Kr$\  is nonempty if and only if \ $\R^\shrp$ \ exists. Therefore, we will implicitly assume
that \ $\R^\shrp$ \ exists.
\end{remark}

\section{{\mathversion{bold}$\Sigma_1$}--Gaps}\label{gaps}

Let \ $M$ \ be transitive model above the reals with a ``cumulative'' hierarchy, say \
$M=\bigcup_{\alpha\in\OR}M_\alpha$ \ where each \ $M_\alpha$ \ is transitive, \
$M_\alpha\subseteq M_\beta$ \ for \ $\alpha\le\beta$, \ and \
$M_\alpha=\bigcup_{\beta<\alpha}M_\beta$ \ for limit \ $\alpha$. \ One can discuss the question of
``when do new \ $\Sigma_1$ \ truths about the reals (that is, about \ $\R\cup\{\R\}$) occur  in \
$M$?'' Suppose, for example, that \ $\varphi$ \ is a \ $\Sigma_1$ \ formula which has the set \
$\R$ \ and an \ $x\in\R$
\ as parameters. If \ $M_\alpha$ \ is such that \ $M_\alpha\models\varphi$ \ and \
$M_\gamma\not\models\varphi$ \ for all \ $\gamma<\alpha$, \ then will say that \ $M_\alpha$ \
witnesses a new \ $\Sigma_1$ \ truth about the reals. Suppose that the ordinals \ $\alpha\le\beta$ \
are such that \ $M_\alpha$ \ and \ $M_{\beta+1}$ \ both witness {\it new\/} \ $\Sigma_1$ \ truths
about the reals. If, in addition, \ $M_\alpha$ \ and \ $M_\beta$ \ both satisfy the {\it same\/} \
$\Sigma_1$ \ truths about the reals, then we will call the interval \ $[\alpha,\beta]$ \ a \
$\Sigma_1(M)$--gap, \ where \ $[\alpha,\beta]$ \ is defined to be the set of ordinals
\ $[\alpha,\beta]=\{\gamma\in\OR : \alpha\le\gamma\le\beta\}$. 

We will now review the concept of a \ $\Sigma_1$--gap in \ $\Lr$. \ Then we
will focus on the notion of a \ $\Sigma_1$--gap in an iterable real premouse.

\subsection{{\mathversion{bold}$\Sigma_1$}--gaps in {\mathversion{bold}$L(\mathbb{R})$}}\label{lrsubsect}

Steel \cite{Steel} develops a fine structure theory for \ $\Lr$ \ and a Levy hierarchy for \
$\Lr$ \ and then solves the problem of finding scales of
minimal complexity in \ $\Lr$. \ Given a set of reals \ $A\in\Lr$, \ using the reflection
properties of the Levy hierarchy for \ $\Lr$, \ Steel identifies the first level \
$\boldface{\Sigma}{n}(J_\alpha(\R))$ \ at which a scale on \ $A$
\ is definable. This level occurs very close to the first ordinal \ $\alpha$ \ such that \ $A\in
J_\alpha(\R)$ \ and for some \ $\Sigma_1$ \ formula (allowing \ $\R$ \ to appear as a
constant) \ $\varphi(v)$ \ and for some \ $x\in\R$, \ one has \ $J_{\alpha+1}(\R)\models \varphi[x]$
\ and yet \ $J_{\alpha}(\R)\not\models \varphi[x]$. \ So in \ $\Lr$ \ there is a close connection
between obtaining scales of minimal complexity and new \ $\Sigma_1$ \ truths about the reals.
Accordingly, Steel introduces the following definition (see \cite[Definition 2.2]{Steel}).

\begin{definition} \label{Lrgap} Let \ $\alpha\le\beta$ \ be ordinals. The interval \
$[\alpha,\beta]$ \ is a \ {\it $\Sigma_1(\Lr)$--gap\/} if and only if  
\be
\item $J_{\alpha}(\R)\prec_1^{\R}J_{\beta}(\R)$
\item $J_{\alpha'}(\R)\nprec_1^{\R}J_{\alpha}(\R)$ \ for all \ $\alpha'<\alpha$
\item $J_{\beta}(\R)\nprec_1^{\R}J_{\beta'}(\R)$ \ for all \ $\beta'>\beta$.
\ee
\end{definition}

\begin{remark} As noted earlier, we always allow the parameter \ $\R$ \ to
appear as a constant in our relevant languages. Thus, the
statement \ $J_{\alpha}(\R)\prec_1^{\R}J_{\beta}(\R)$ \ is in fact equivalent to the statement
\ $J_{\alpha}(\R)\prec_1^{\R\, \cup \{\R\}}J_{\beta}(\R)$.
\end{remark}

Let \ $\delta\in\OR$ \ be the least such that \ $J_{\delta}(\R)\prec_1^{\R}\Lr$. \ One can show
that the \ $\Sigma_1(\Lr)$--gaps partition \ $\delta$. \ 
Moreover, since  \ $J_{\delta}(\R)\prec_1^{\R}J_{\beta}(\R)$ \ for all
ordinals \ $\beta>\delta$, \ it follows that \ $\delta$ \  starts a \
$\Sigma_1(\Lr)$--gap which has ``no end''. In Section~\ref{minimal}, however, we will show that this particular gap will have a ``proper'' ending in any iterable real premouse. 
\subsection{{\mathversion{bold}$\Sigma_1$}--gaps in {\mathversion{bold}$K(\mathbb{R})$}}
Since \ $\Kr$ \ is the union of real 1--mice, the notion of a \ $\Sigma_1$--gap in \ $\Kr$ \ reduces to discussing such gaps
in iterable real premice. We recall that an iterable real premouse \ $\nouse=(N,\R,\kappa,\mu)$ \ is a structure,
consisting of sets constructible above the reals relative to the measure \ $\mu$, \ with a Jensen hierarchy similar to that
in subsection~\ref{jensenL(R)} (see Section 1 of \cite{Crcm}). We denote the \ $\alpha^{\,\ul{\text{th}}}$ \ level of this
hierarchy by \ $\nouse^\alpha$. \ This avoids any confusion with the notation \ $\nouse_\alpha$ \ which denotes the
\ $\alpha^{\,\ul{\text{th}}}$ \ premouse iterate of \ $\nouse$.
\ In Section~\ref{minimal} we shall identify the pointclasses of the form \ $\boldface{\Sigma}{n}(\nouse^\alpha)$ \ which have
the scale property. 
\subsubsection*{{\mathversion{bold}$\Sigma_1$}--gaps in iterable real premice}
Let \ $\nouse=(N,\R,\kappa,\mu)$ \ be an iterable real premouse and let \ $\mouse$ \ be the
iterable real premouse defined by \ $\mouse=(M,\R,\kappa,\mu)=\nouse^{\kappa+1}$. \ It can be shown
that every set of reals in \ $\Lr$ \ is in \ $M$ \ and that for all \ $\alpha\le\kappa$, \
$J_\alpha(\R)=M^\alpha$. \ In fact, letting \ $\delta$ \ be as in the above subsection~\ref{lrsubsect}, one can show that \
$\delta<\kappa$ \ and \ $J_{\delta}(\R)\nprec_1^{\R}\mouse$. \ Thus, in \ $\mouse$ \ the
ordinal \ $\kappa$ \ is the end of the \ $\Sigma_1$--gap that began with \ $\delta$ \ and  the
ordinal \ $\kappa+1$ \ starts a new \ $\Sigma_1$--gap.  Furthermore, one can show that there is a
set of reals \ $A\notin\Lr$ \ such that \ $A$ \ has a \ $\Sigma_1(\mouse)$ \
scale.\footnote{Namely, \ $A=\R^\shrp$. The iterable real premouse $\mouse$ is
``$\R$--sharplike'' (see footnote~\ref{sharpnote}).}  Consequently, the connection between new scales
and  \ $\Sigma_1$ \ truths continues.

We will show that in iterable real premice there is a close connection between obtaining scales of minimal complexity and new \
$\Sigma_1$ \ truths about the reals. The following definition is a straightforward generalization of Definition
\ref{Lrgap}. 
\begin{definition} \label{gap} Let \ $\nouse=(N,\R,\kappa,\mu)$ \ be an iterable real premouse. Let
\
$\kappa<\alpha\le\beta\le\widehat{\OR}^{\,\nouse}$. \ We shall say that
the interval \ $[\alpha,\beta]$ \ is a \ {\it $\Sigma_1(\nouse)$--gap\/} \ if and only if
\be
\item $\nouse^{\alpha}\prec_1^{\R}\nouse^{\beta}$
\item $\nouse^{\alpha'}\nprec_1^{\R}\nouse^{\alpha}$ \ for all \ $\alpha'<\alpha$
\item $\nouse^{\beta}\nprec_1^{\R}\nouse^{\beta'}$  \ for all \ $\beta'>\beta$ \ where \ $\beta'\le\widehat{\OR}^{\,\nouse}$.
\ee
\end{definition}

When \ $[\alpha,\beta]$ \ is a \ $\Sigma_1(\nouse)$--gap we shall say that \ $\alpha$ \
{\it begins\/} \ the gap and that \ $\beta$ \ {\it ends\/} \ the gap. In addition, if \
$\alpha<\widehat{\OR}^{\,\nouse}$ \ then we shall say that \ $\alpha$ \ {\it properly
begins\/} the gap and, when the context is clear, we will say that \ $\alpha$ \ is proper. When \
$\beta<\widehat{\OR}^{\,\nouse}$
\ we will say that \ $\beta$ \ {\it properly ends\/} the gap and, when the context is clear, we shall say that \ $\beta$ \
is proper. 
\begin{remark} Two observations concerning Definition~\ref{gap}: 
\be
\item We allow for the possibility that 
\ $[\alpha,\alpha]$ \ is a \ $\Sigma_1(\nouse)$--gap \ when
\ $\alpha\le\widehat{\OR}^{\,\nouse}$. 
\item The ordinal \ $\delta=\widehat{\OR}^{\,\nouse}$ \ always ends a $\Sigma_1(\nouse)$--gap. More
specifically, either (i) $[\delta,\delta]$ \ is a \ $\Sigma_1(\nouse)$--gap, or (ii) $[\alpha,\delta]$ \ is a \
$\Sigma_1(\nouse)$--gap for some \ $\alpha<\delta$.
\ee
\end{remark}
\begin{remark}\label{why} 
Suppose that an iterable real premouse \ $\mouse$ \ is a proper extension of \ $\nouse$. \ If \ $\alpha$ \
begins a $\Sigma_1(\nouse)$--gap then \ $\alpha$ \ will also begin a  $\Sigma_1(\mouse)$--gap.
Also, if \ $[\alpha,\beta]$ \ is a \ $\Sigma_1(\nouse)$--gap  and \
$\beta<\widehat{\OR}^{\,\nouse}$, \ then \ $[\alpha,\beta]$ \ will likewise be a
\ $\Sigma_1(\mouse)$--gap.  If \ $\beta=\widehat{\OR}^{\,\nouse}$, \ then
it is possible for the interval \ $[\alpha,\beta]$ \ to fail to be a \ $\Sigma_1(\mouse)$--gap,
because the end of the corresponding \ $\Sigma_1(\mouse)$--gap  may be greater than \ $\beta$.
\end{remark}

\begin{definition}\label{internal} 
Let \ $\nouse=(N,\R,\kappa,\mu)$ \ be an iterable real premouse.  For \
$\gamma<\widehat{\OR}^{\,\nouse}$ \ we shall say that \ $\boldface{\Sigma}{n}(\nouse^\gamma)$ \ and
\ $\boldface{\Pi}{n}(\nouse^\gamma)$, \ for \ $n\ge 1$, \ are {\it internal\/} pointclasses.
On the other hand, we shall call \ $\boldface{\Sigma}{n}(\nouse)$ \ and \
$\boldface{\Pi}{n}(\nouse)$, \ for \ $n\ge 1$, \  {\it external\/} pointclasses.
\end{definition}

\section{Complexity of Scales in an Iterable Real Premouse}\label{minimal}
In this section we shall assume that \ $\nouse=(N,\R,\kappa,\mu)$ \ is an
iterable real premouse. We shall use the notions developed in the previous sections to identify
those levels of the Levy hierarchy for \ $\nouse$ \ which have the scale property.

\subsection{Scales at the beginning of a gap}

Recall that an iterable real premouse \ $\mouse$ \ is a 1--mouse if \
$\omega\rho_{\mouse}\le \kappa^{\mouse}$ \ (see Definition~\ref{onemouse}). 
Definition~{\hulldef} of \cite{Part1} describes the $\Sigma_n$ hull of a premouse. Let \ $\mouse$ \  be a 1--mouse, let \ 
$\mathcal{H} = \Hull_1^{\mouse}({\R\cup \omega\rho_{\mouse}} \cup \{p_{\mouse}\})$, \ and let  \ $\core$ \  be the transitive
collapse of \ $\mathcal{H}$. \ It follows that \ $\core$ \ is a real premouse, denoted by  \ $\core(\mouse)$.

\begin{theorem}\label{beggap} Suppose that \
$\alpha$ \ begins a \ $\Sigma_1(\nouse)$--gap. Let \ $\C=(C,\R,\kappa^\C,\mu^\C)$ \ be the
transitive collapse of \ $\Hull_1^{\,\nouse^\alpha}(\R)$. \ Then 
\newcounter{stuff}
\newcounter{stuffit}
\begin{list}{{\upshape (\roman{stuff})}}
{\setlength{\labelwidth}{20pt}
\usecounter{stuff}}
\item $\C$ \ is an iterable real premouse and \ $\C_\xi=\nouse^\alpha$ \ for some ordinal \
$\xi$.
\item There exists a set of reals \ $A$ \ such that \ $A\in\Sigma_1(\nouse^\alpha)$ \ and \
$A\notin N^\alpha$. 
\item $\nouse^\alpha$ \ is a real 1--mouse, \ $\core(\nouse^\alpha)=\C$, \
$\rho^{}_{{\nouse^\alpha}}=1$ \ and \ $p^{}_{{\nouse^\alpha}}=\emptyset$.
\item \setcounter{stuffit}{\value{stuff}}
$\boldface{\Sigma}{1}(\nouse^\alpha)=\Sigma_1(\nouse^\alpha,\R)$,
\ as pointclasses.
\item $\pow(\R)\cap N^\alpha\subset\Sigma_1(\nouse^\alpha,\R)$.
\end{list}
\end{theorem}

\begin{proof} Let \ $\alpha$, \ $\nouse$ \ and \ $\C$ \ be as stated in the theorem. We prove items (i)--(iv).

(i) Let \ $\sigma\colon\C\mapsigma{1}\Hull_1^{\,\nouse^\alpha}(\R)$ \ be the inverse of the
collapse map \ $\pi\colon \Hull_1^{\,\nouse^\alpha}(\R)\to\C$. \ It follows that \ $\C$ \ is a real premouse. Because \ $\sigma$ \ is a \
$\Sigma_1$
\ embedding, we have that the induced map \
$\sigma\colon F^{\C}\to F^{\nouse^\alpha}$ \ is \ $\le$-extendible (see Definition~{\defttho} of \cite{Part1}).
Thus, Theorem~{\thmttho} of \cite{Part1} implies that \  $\C$ \ is an iterable real premouse. Note that \
$\C=\Hull_1^{\,\C}(\R)$. \ We now conclude from Theorem 2.32 of \cite{Crcm} that there
is an ordinal \ $\xi$ \ such that the \ $\xi^{\,\ul{\text{th}}}$ \ premouse iterate \ $\C_\xi$ \ is an initial segment
of \ $\nouse^\alpha$. \  Because \ $\pi_{0\xi}\colon\C\mapsigma{1}\C_\xi$, \ where \ $\pi_{0\xi}$ \
is the  premouse embedding of \ $C$ \ into \ $\C_\xi$, \ it follows that \
$\C_\xi\prec_1^\R\nouse^\alpha$ \ and, \ since \ $\alpha$ \ begins a \
$\Sigma_1(\nouse)$--gap, we see that \ $\C_\xi=\nouse^\alpha$.

(ii) We first define a set of reals \ $A$ \ such that \
$A\in\Sigma_1(\C)$ \ and \ $A\notin C$. \ Let \ $g=h_\C\restriction(\R\times\R)$ \ be the partial map obtained by restricting
the canonical \ $\Sigma_1(\C)$ \ Skolem function \ $h_\core$ \ (in \cite{Part1} see Definition~{\skolemdef} and Lemma
{\skolemimage}).
\ Because
\
$\C=\Hull_1^{\,\C}(\R)$, \ it follows that
\ $g\colon\R\times\R\maps{onto} C$. \ Define \ $A\subseteq \R$ \ by
\[x\in A \iff x\notin g(x_0,x_1).\]
Note that \ $A\in\Sigma_1(\C)$ \ and, since \ $g$ \ is onto,  \ $A\notin C$. \ 
Because \ $A\in\Sigma_1(\C)$, \ we have that \ $A\in\Sigma_1(\nouse^\alpha)$. \ Also, since \
$\C_\xi=\nouse^\alpha$, \ Lemma 2.11(4) of \cite{Crcm} implies that \ $A\notin
N^\alpha$. \

(iii) It  follows from (i) and (ii) that \ $\nouse^\alpha$ \ is a real 1--mouse, \
$\C=\core(\nouse^\alpha)$, \ $\rho^{}_{{\nouse^\alpha}}=1$ \ and \
$p^{}_{{\nouse^\alpha}}=\emptyset$.

(iv) Because \ $\sigma\colon\C\mapsigma{1}\Hull_1^{\,\nouse^\alpha}(\R)$, \ we see that 
\ $\Sigma_1(\nouse^\alpha,\R)=\Sigma_1(\C,\R)$. \ In addition, since \
$\C=\Hull_1^{\,\C}(\R)$, \ we have that \ $\boldface{\Sigma}{1}(\C)=\Sigma_1(\C,\R)$.
Since \ $\C_\xi=\nouse^\alpha$, \ Corollary 2.14 of \cite{Crcm} implies that \
$\boldface{\Sigma}{1}(\nouse^\alpha)=\boldface{\Sigma}{1}(\C)$. \ Therefore, \ 
$\boldface{\Sigma}{1}(\nouse^\alpha)=\Sigma_1(\nouse^\alpha,\R)$.

(v) This follows immediately from (ii) and (iv).
\end{proof}

\begin{theorem}\label{alphascales}  Suppose that \
$\alpha$ \ begins a \ $\Sigma_1(\nouse)$--gap. If \ $\nouse^\alpha \models
\AD$, \ then \ $\boldface{\Sigma}{1}(\nouse^\alpha)$ \ has the scale property. 
\end{theorem}

\begin{proof} Corollary~\ref{relscales} states that \ $\Sigma_1(\nouse^\alpha,\R)$ \ has the scale
property. Therefore, \ $\boldface{\Sigma}{1}(\nouse^\alpha)$ \ has the scale property by (iv) of
Theorem~\ref{beggap}.\end{proof}

When \ $\alpha$ \ begins a \ $\Sigma_1(\nouse)$--gap, we conclude that \ $\boldface{\Pi}{1}(\nouse^\alpha)$ \ does not have
the scale property (see \cite[4B.13]{Mosch}). The next classes to consider are the pointclasses \
$\boldface{\Sigma}{n}(\nouse^\alpha)$ \ where \
$n>1$ \ and \ $\alpha$ \ properly begins a \ $\Sigma_1(\nouse)$--gap.

\begin{definition}\label{bool} Let \ $\Gamma\subseteq\pow(\R)$ \ be a pointclass. Define
\[\mathfrak{B}(\Gamma)=\{ A\subseteq\R : \text{$A$ is a boolean combination of $\Gamma$
sets}\}.\]
\end{definition}
$\mathfrak{B}(\Gamma)$ \ is a pointclass that is closed under the operations of intersection,
union and complement.
\begin{definition}\label{realquantifier} Given a  pointclass \ $\Gamma$, \ define \ $\Sigma_n^*(\Gamma)$
\ by the following induction on \ $n\in\omega$,
\begin{align*}
\Sigma_0^*(\Gamma)&=\Gamma\\
\Pi_n^*(\Gamma)&=\lnot\Sigma_n^*(\Gamma)\\
\Sigma_{n+1}^*(\Gamma)&=\exists^{\R}\Pi_n^*(\Gamma)
\end{align*}
\end{definition}

\begin{remark}\label{spthm:rmk} Let \ $\mouse=(M,\R,\kappa,\mu)$ \ be a real premouse. Applying
Definition
\ref{realquantifier} to the pointclass \ $\Gamma=\Sigma_1(\mouse,\R)$ \ we obtain the following
equations:
\begin{align*}
\Sigma_0^*(\Sigma_1(\mouse,\R))&=\Sigma_1(\mouse,\R))\\
\Pi_0^*(\Sigma_1(\mouse,\R))&=\Pi_1(\mouse,\R))\\
\Sigma_{1}^*(\Sigma_1(\mouse,\R))&=\exists^{\R}\Pi_1(\mouse,\R)\\
\Pi_{1}^*(\Sigma_1(\mouse,\R))&=\forall^{\R}\Sigma_1(\mouse,\R)\\
\Sigma_{2}^*(\Sigma_1(\mouse,\R))&=\exists^{\R}\forall^{\R}\Sigma_1(\mouse,\R)\\
\Pi_{2}^*(\Sigma_1(\mouse,\R))&=\forall^{\R}\exists^{\R}\Pi_1(\mouse,\R).
\end{align*}
Consequently, assuming enough determinacy, Corollary~\ref{relscales} and the
Second Periodicity Theorem of Moschovakis \cite[Theorems 6C.2 and  6C.3]{Mosch} imply that the pointclasses
\[\Sigma_0^*(\Sigma_1(\mouse,\R)), \
\Pi_{1}^*(\Sigma_1(\mouse,\R)), \ 
\Sigma_{2}^*(\Sigma_1(\mouse,\R))\] have the scale property.
\end{remark}
We shall motivate our next proposition with two examples. Let 
\ $\mouse=(M,\R,\kappa,\mu)$ \ be a real premouse. Suppose that \ $\mouse=\Hull_1^{\,\mouse}(\R)$ \
and let \ $g\colon\R\to M$ \ be the {\sl partial\/} function defined by \ $g(x)=h_\mouse(x_0,x_1)$.
\ Thus, \ $g$ \ is onto \ $M$ \ with \ $\Sigma_1(\mouse)$ \ graph.  Let \ $D=\dom(g)$ \ and note
that \ $D\in\Sigma_1(\mouse)$.

\begin{example}\label{exmone} Now consider the pointclass
\ $\boldface{\Sigma}{2}(\mouse)$. \ Let \ $A\in\boldface{\Sigma}{2}(\mouse)$. \ Then there is a
\ $\Pi_1$
\ formula \ $\varphi(u,v,w)$ \ and a \ $k\in M$ \ such that 
\[x\in A \iff \mouse\models (\exists u)\varphi(u,x,k)\]
for all \ $x\in\R$. \ 
Let \ $y\in\R$ \ be such that \ $g(y)=k$, \ where \ $g$ \ is defined above.
 \ Since
\begin{alignat*}{1}
x\in A &\iff \mouse\models (\exists u)\varphi(u,x,g(y))\\
&\iff \mouse\models (\exists z\in\R)[z\in D\land
(\forall a)(\forall b)((a=g(z)\land b=g(y))
\rightarrow
\varphi(a,x,b))],
\end{alignat*}
we conclude that \ $A\in \exists^{\R}(\mathfrak{B}(\Sigma_1(\mouse,\R)))$. \ By
definition,
\[\exists^{\R}(\mathfrak{B}(\Sigma_1(\mouse,\R)))=\Sigma_1^*(\mathfrak{B}(\Sigma_1(\mouse,\R)))\]
and so, \ $\boldface{\Sigma}{2}(\mouse)\subseteq
\Sigma_{1}^*(\mathfrak{B}(\Sigma_1(\mouse,\R)))$.
\end{example}
\begin{example}\label{exmtwo} In addition, consider the pointclass \ $\boldface{\Sigma}{3}(\mouse)$. \ Let \
$A\in\boldface{\Sigma}{3}(\mouse)$. \ Then there is a \ $\Sigma_1$ \
formula \ $\psi(t,u,v,w)$ \ and a \ $k\in M$ \ such that 
\[x\in A \iff \mouse\models (\exists t)(\forall u)\psi(t,u,x,k)\]
for all \ $x\in\R$. \ 
Let \ $y\in\R$ \ be such that \ $g(y)=k$. \ Since
\begin{alignat*}{1}
x\in A&\iff \mouse\models (\exists t)(\forall u)\psi(t,u,x,g(y))\\
&\iff \mouse\models (\exists p\in\R)(\forall z\in\R)[p\in D\land (z\in D \rightarrow \\ 
&\qquad \qquad (\exists a)(\exists b)(\exists c)(a=g(p)\land b=g(z) \land c=g(y)\land
\psi(a,b,x,c)))]
\end{alignat*}
we conclude that \ $A\in
\exists^{\R}\forall^{\R}(\mathfrak{B}(\Sigma_1(\mouse,\R)))$. \ By definition,  
\[\exists^{\R}\forall^{\R}(\mathfrak{B}(\Sigma_1(\mouse,\R)))=
\Sigma_2^*(\mathfrak{B}(\Sigma_1(\mouse,\R)))\]
and so, \ $\boldface{\Sigma}{3}(\mouse)\subseteq
\Sigma_{2}^*(\mathfrak{B}(\Sigma_1(\mouse,\R)))$. 
\end{example}
This completes our two examples. The following proposition should now be clear. 

\begin{proposition}\label{propp} Let \ $\mouse=(M,\R,\kappa,\mu)$ \ be a real premouse. Assume that \
$\mouse=\Hull_1^{\,\mouse}(\R)$. \ Then \
$\boldface{\Sigma}{n+1}(\mouse)\subseteq
\Sigma_{n}^*(\mathfrak{B}(\Sigma_1(\mouse,\R)))$ \ for all \ $n\ge 1$.
\end{proposition}

\begin{lemma}\label{subsetrel}  Let \ $\nouse=(N,\R,\kappa,\mu)$ \ be an iterable real premouse and
suppose that \ $\alpha$ \ begins a \ $\Sigma_1(\nouse)$--gap. Then 
\[\boldface{\Sigma}{n+1}(\nouse^\alpha)\subseteq
\Sigma_{n}^*(\mathfrak{B}(\Sigma_1(\nouse^\alpha,\R)))\] for all \ $n\ge 1$.
\end{lemma}
\begin{proof} Since \ $\alpha$ \ begins a \ $\Sigma_1(\nouse)$--gap. Let \ $\C=(C,\R,\kappa^\C,\mu^\C)$ \
be the transitive collapse of \ $\H=\Hull_1^{\,\nouse^\alpha}(\R)$.  
\ By Theorem~\ref{beggap}, there is a premouse iterate \ $\C_\xi$ \ such that \ $\C_\xi=\nouse^\alpha$. \ Let \ $n\ge 1$. \ 
\begin{claim} $\Sigma_1(\C,\R)=\Sigma_1(\nouse^\alpha,\R)$ \
and \ $\boldface{\Sigma}{n+1}(\C)=\boldface{\Sigma}{n+1}(\nouse^\alpha)$
\ as pointclasses.
\end{claim}
\begin{proof}[Proof of Claim] Clearly, \ $\Sigma_1(\C,\R)=\Sigma_1(\nouse^\alpha,\R)$.
Corollary 2.20 of \cite{Crcm} implies that \
$\boldface{\Sigma}{n+1}(\nouse^\alpha)\subseteq
\boldface{\Sigma}{n+1}(\C)$. \ Let \ $\H^\alpha=\Hull_1^{\,\nouse^\alpha}(\R)$. \ Note that \ $\C$ \ and \ $\H^\alpha$ \ are
isomorphic structures and thus, \ $\boldface{\Sigma}{n+1}(\C)=\boldface{\Sigma}{n+1}(\H^\alpha)$. \ Because \ $\H^\alpha$ \ is
\ $\Sigma_1$ \ definable over \ $\nouse^\alpha$, \ we have that \
$\boldface{\Sigma}{n+1}(\H^\alpha)\subseteq\boldface{\Sigma}{n+1}(\nouse^\alpha)$. \ Therefore, \
$\boldface{\Sigma}{n+1}(\C)=\boldface{\Sigma}{n+1}(\nouse^\alpha)$.
\end{proof}
Proposition~\ref{propp} asserts that \ $\boldface{\Sigma}{n+1}(\C)\subseteq
\Sigma_{n}^*(\mathfrak{B}(\Sigma_1(\C,\R)))$. \ Thus, the Claim implies that \
$\boldface{\Sigma}{n+1}(\nouse^\alpha)\subseteq
\Sigma_{n}^*(\mathfrak{B}(\Sigma_1(\nouse^\alpha,\R)))$. 
\end{proof}
\begin{definition} Let \ $\mouse$ \ be a real premouse. We say that \ $\mouse$ \ is \
{\it $\R$--collectible\/} if for every \ $\Sigma_0$ \ formula \ $\varphi(x,v)$
\ (allowing arbitrary parameters in $M$) we have that
\[\mouse\models (\forall x\in\R)(\exists v)\varphi(x,v)\text{ \ implies \ }\mouse\models
(\exists w)(\forall x\in\R)(\exists v\in w)\varphi(x,v).\]
\end{definition}

\begin{lemma}\label{collection.equiv} Let \ $\nouse=(N,\R,\kappa,\mu)$ \ be an iterable real
premouse and suppose that \ $\alpha$ \ begins a \ $\Sigma_1(\nouse)$--gap. Then
\ $\nouse^\alpha$ is \ $\R$--collectible if and only if \ $\core(\nouse^\alpha)$ is \
$\R$--collectible.
\end{lemma}
\begin{proof}
Suppose that \ $\alpha$ \ begins a \ $\Sigma_1(\nouse)$--gap. Let \ $\C=(C,\R,\kappa^\C,\mu^\C)$ \
be the transitive collapse of \ $\H=\Hull_1^{\,\nouse^\alpha}(\R)=(H,\R,\kappa^\H,\mu^\H)$. \ By
Theorem~\ref{beggap}, \
$\C=\core(\nouse^\alpha)$ \ and \ $\C_\xi=\nouse^\alpha$ \ for some ordinal \
$\xi$.

\smallskip
$(\Longrightarrow)$. \ Assume \ $\nouse^\alpha$ is \ $\R$--collectible. We will show that \ $\H$ \ is
\ $\R$--collectible. It will then follow that \ $\core(\nouse^\alpha)$ is \ $\R$--collectible.
Note that  \ $\H\prec_1\nouse^\alpha$ \ and \ $\R\cup\{\R\}\in H$. \ For simplicity let \
$\varphi(x,v_1,v_2)$ \ be \ $\Sigma_0$ \ and let \ $a\in H$. \ Suppose that
\begin{equation}\H\models (\forall x\in\R)(\exists
v)\varphi(x,v,a).\label{assump0}\end{equation}
Since \
$\H\prec_1\nouse^\alpha$, \ (\ref{assump0}) implies that
\[\nouse^\alpha\models (\forall x\in\R)(\exists v)\varphi(x,v,a).\]
Then, because \ $\nouse^\alpha$ \ is \ $\R$--collectible,
\[\nouse^\alpha\models(\exists w)(\forall x\in\R)(\exists v\in w)\varphi(x,v,a).\]
Since \ $\H\prec_1\nouse^\alpha$, \ it  follows that
\[\H\models(\exists w)(\forall x\in\R)(\exists v\in w)\varphi(x,v,a).\]
Therefore,  \ $\core(\nouse^\alpha)$ is \ $\R$--collectible.

\smallskip
$(\Longleftarrow)$. \ Assume \ $\C=\core(\nouse^\alpha)$ is \ $\R$--collectible. Let 
\[\premouseiteration{\C}{\gamma}{\beta}\]
be the premouse iteration of \ $\C$. \ Note that \ $\R\cup\{\R\}\in C_\gamma$ \ for all \
$\gamma\in\OR$ \ (see  \cite[Lemma 2.11(3)]{Crcm}). We will show by induction on
\ $\lambda$, \ that \ $\C_\lambda$ \ is \ $\R$--collectible. We can then conclude that \ $\nouse^\alpha$
\ is \ $\R$--collectible, since \ $\C_\xi=\nouse^\alpha$. \ Now, for \ $\lambda=0$, \ $\C_0=\C$. \
Hence, \ $\C_0$ \ is $\R$--collectible by assumption. 

{\sc Successor Case:} Let \ $\lambda\ge 0$ \ and assume that \ $\C_\lambda$ \ is \ $\R$--collectible.
We will show that \ $\C_{\lambda+1}$ \ is \ $\R$--collectible. Let \ $\pi=\pi_{\lambda,\lambda+1}$. \
Thus, \ $\pi\colon\C_\lambda \updownmap{cofinal}{1}\C_{\lambda+1}$. \ For simplicity let \
$\varphi(x,v_1,v_2)$ \ be \ $\Sigma_0$ \ and let \ $a\in C_{\lambda+1}$. \ Now suppose that
\begin{equation}\C_{\lambda+1}\models (\forall x\in\R)(\exists
v)\varphi(x,v,a).\label{assump}\end{equation} Since \ $\pi$ \ is cofinal, \ (\ref{assump}) implies
\[(\forall x\in\R)(\exists b\in C_\lambda)\left[\,\C_{\lambda+1}\models (\exists
v\in\pi(b))\varphi(x,v,a)\,\right].\]
By Lemma 2.8(3) of \cite{Crcm}, there is a function \ $f\in C_\lambda$ \ such that \
$\pi(f)(\kappa_\lambda)=a$. \ Thus by Theorem 2.4 of \cite{Crcm}, we conclude that
\[(\forall x\in\R)(\exists b_1,b_2\in C_\lambda)\left[\,\C_{\lambda}\models
\!\left(b_0=\{\eta\in\underline{\kappa} : (\exists v\in b_1)\varphi(x,v,f(\eta))\}\land
b_0\in\mu\right)\,\right].\]
By our induction hypothesis, there exists a \ $w\in C_\lambda$ \ such that
\[\C_{\lambda}\models(\forall x\in\R)(\exists b_1,b_2\in w)\left[\
b_0=\{\eta\in\underline{\kappa} : (\exists v\in b_1)\varphi(x,v,f(\eta))\}\land
b_0\in\mu\,\right].\]
It follows that 
\[\C_{\lambda+1}\models (\forall x\in\R)(\exists
v\in\pi(w))\varphi(x,v,a).\]
This argument shows that \ $\C_{\lambda+1}$ \ is \ $\R$--collectible.

{\sc Limit Case:} Let \ $\lambda\ge 0$ \ be a limit ordinal and assume that \ $\C_\nu$ \ is \
$\R$--collectible for all \ $\nu<\lambda$. \ We will show that \ $\C_{\lambda}$ \ is \
$\R$--collectible. For simplicity let \
$\varphi(x,v_1,v_2)$ \ be \ $\Sigma_0$ \ and let \ $a\in C_{\lambda}$. \ Suppose that
\begin{equation}\C_{\lambda}\models (\forall x\in\R)(\exists
v)\varphi(x,v,a).\label{assump2}\end{equation} 
Because \ $\C_\lambda$ \ is a direct limit, there is a an ordinal \ $\nu<\lambda$ \ and an
element \ $\hat{a}\in\C_\nu$ \ such that \
$\pi_{\nu\lambda}(\hat{a})=a$. \ Since \
$\pi_{\nu\lambda}\colon\C_\nu\updownmap{}{1}\C_\lambda$, \ we see from
(\ref{assump2}) that
\[\C_\nu\models (\forall x\in\R)(\exists
v)\varphi(x,v,\hat{a}).\]
By our induction hypothesis, there exists a \ $w\in C_\nu$  \ such that
\[\C_\nu\models (\forall x\in\R)(\exists
v\in w)\varphi(x,v,\hat{a}).\]
Therefore, 
\[\C_\lambda\models (\forall x\in\R)(\exists
v\in\pi_{\nu\lambda}(w))\varphi(x,v,a).\]
This argument shows that \ $\C_\lambda$ \ is \ $\R$--collectible.
\end{proof}
\begin{corollary}\label{total} Let \ $\nouse=(N,\R,\kappa,\mu)$ \ be an iterable real
premouse and suppose that \ $\alpha$ \ begins a \ $\Sigma_1(\nouse)$--gap. If \ $\nouse^\alpha$ is not \ $\R$--collectible,
then there is a total function  \ $k\colon\R\maps{cofinal}\OR^{\nouse^\alpha}$ \  whose graph is
\ $\Sigma_1(\nouse^\alpha,\R)$.
\end{corollary}
\begin{proof} Assume that \ $\nouse^\alpha$ is not \ $\R$--collectible. Let \ $\C=\core(\nouse^\alpha)$ \ and let \ $\xi$ \ be
such that \ $\C_\xi=\nouse^\alpha$ \ (see (i) of Theorem~\ref{beggap}). By Lemma~\ref{collection.equiv}, \ $\C$ \ is not \
$\R$--collectible. It follows that there is a total function \ $k'\colon\R\maps{cofinal}\OR^\C$ \  whose graph is \
$\boldface{\Sigma}{1}(\C)$. \ Since there is a partial \ $\Sigma_1(C)$ \ map of \ $\R$ \ onto \
$C$, \ we see that the graph of \ $k'$ \ is \  $\Sigma_1(\C,\R)$. \ Because \
$\pi_{0\xi}\colon\C\updownmap{cofinal}{1}\nouse^\alpha$, \ where \ $\pi_{0\xi}$ \ is the  premouse
embedding of \ $C$ \ into \ $\C_\xi=\nouse^\alpha$, \ we conclude that there is a total function  \
$k\colon\R\maps{cofinal}\OR^{\nouse^\alpha}$ \  whose graph is \ $\Sigma_1(\nouse^\alpha,\R)$ \
($k$ \ is just the interpretation of \ $k'$ \ in \ $\nouse^\alpha$).
\end{proof}

\begin{definition} Let \ $\nouse=(N,\R,\kappa,\mu)$ \ be an iterable real premouse and suppose that
\ $\alpha$ \ begins a \ $\Sigma_1(\nouse)$--gap. Then  we shall say that
\begin{itemize}
\item $\alpha$ \ is {\it collectible\/}  if and only if  \ $\nouse^\alpha$ \ is \
$\R$--collectible, 
\item $\alpha$ \ is {\it uncollectible\/} if and only if  \ $\nouse^\alpha$ \ is not \
$\R$--collectible.
\end{itemize}
\end{definition}

\begin{rmk} Suppose that \ $\nouse$ \ is an iterable real premouse and \ $[\alpha, \beta]$ \ is a \ $\Sigma_1(\nouse)$--gap. If
\ $\alpha<\beta$ \ then \ $\alpha$ \ is collectible; however, if \ $\alpha=\beta$ \ then \ $\alpha$ \ can be  uncollectible.
\end{rmk}

\begin{lemma}\label{pntclsequal}  Suppose that \ $\alpha$ \ begins a \ $\Sigma_1(\nouse)$--gap. If \
$\alpha$
\ is uncollectible, then 
\begin{alignat*}{1}
\boldface{\Sigma}{n+1}(\nouse^\alpha) &=
\Sigma_{n}^*(\Sigma_1(\nouse^\alpha,\R))\tag{a}\\
\boldface{\Pi}{n+1}(\nouse^\alpha) &= \Pi_{n}^*(\Sigma_1(\nouse^\alpha,\R))\tag{b}
\end{alignat*}
(as pointclasses) for all \ $n\ge 0$.
\end{lemma}

\begin{proof}
By (\roman{stuffit}) of Theorem~\ref{beggap}, \
$\boldface{\Sigma}{1}(\nouse^\alpha)=\Sigma_1(\nouse^\alpha,\R)$. \ Thus, the Lemma holds for \
$n=0$, \ by Definition~\ref{realquantifier}. So we assume \ $n\ge 1$.

It should be clear that
\begin{alignat*}{1}
\boldface{\Sigma}{n+1}(\nouse^\alpha) &\supseteq
\Sigma_{n}^*(\Sigma_1(\nouse^\alpha,\R))\\
\boldface{\Pi}{n+1}(\nouse^\alpha) &\supseteq \Pi_{n}^*(\Sigma_1(\nouse^\alpha,\R))
\end{alignat*}
and thus we shall show that
\begin{alignat*}{1}
\boldface{\Sigma}{n+1}(\nouse^\alpha) &\subseteq \Sigma_{n}^*(\Sigma_1(\nouse^\alpha,\R))\\
\boldface{\Pi}{n+1}(\nouse^\alpha) &\subseteq \Pi_{n}^*(\Sigma_1(\nouse^\alpha,\R)).
\end{alignat*}
In fact, it is sufficient to prove that \ $\boldface{\Sigma}{n+1}(\nouse^\alpha) \subseteq
\Sigma_{n}^*(\Sigma_1(\nouse^\alpha,\R))$.
\ Now, to simplify the notation slightly, for any real premouse \ $\mouse$ \ let
\[\mathfrak{B}_1(\mouse)=
\mathfrak{B}(\Sigma_1(\mouse,\R)).\]
By Lemma~\ref{subsetrel}, \ 
$\boldface{\Sigma}{n+1}(\nouse^\alpha) \subseteq
\Sigma_{n}^*(\mathfrak{B}_1(\nouse^\alpha))$.

\begin{claim} $\mathfrak{B}_1(\nouse^\alpha)\subseteq \forall^\R\Sigma_1(\nouse^\alpha,\R)\cap
\exists^\R\Pi_1(\nouse^\alpha,\R)$.
\end{claim}

\begin{proof}[Proof of Claim] Because \ $\mathfrak{B}_1(\nouse^\alpha)$ \ is closed under
complementation, it is sufficient to show that
\begin{equation}\mathfrak{B}_1(\nouse^\alpha)\subseteq\exists^\R\Pi_1(\nouse^\alpha,\R).
\tag{$*$}\end{equation}
To show ($*$),  it is enough to establish that \
$\Sigma_1(\nouse^\alpha,\R)\subseteq\exists^{\R}\Pi_1(\nouse^\alpha,\R)$. \ Let \
$Q\in\Sigma_1(\nouse^\alpha,\R)$, \ say 
\[y\in Q \iff \nouse^\alpha\models \varphi(y,z)\]
for all \ $y\in\R$, \ where \ $\varphi\in\Sigma_1$ \ and \ $z\in\R$.  \ We will show that \
$Q\in\exists^{\R}\Pi_1(\nouse^\alpha,\R)$. \ Corollary
\ref{total} implies that there is a total function  \ $k\colon\R\maps{cofinal}\OR^{\nouse^\alpha}$ \  whose
graph is \ $\Sigma_1(\nouse^\alpha,\R)$ \ and such that \ $k(x)>\kappa$ \ for all \ $x\in\R$.  

Let \ $\langle S_\gamma\ : \ \gamma <
\OR^{\nouse^\alpha} \rangle$ \ be the \ $\Sigma_1(\nouse^\alpha)$ \ increasing sequence
of transitive sets (see Section 1 of \cite{Crcm}) such that 
\bi
\item $\R\in S_\gamma\in N^\alpha$, \ for all \ $\eta \in \OR^{\nouse^\alpha}$
\item $N^\alpha=\bigcup_{\eta \in \zeta}S_\gamma$, \ where \
$\zeta=\OR^{\nouse^\alpha}$.
\ei
For all ordinals \ $\eta$ \ such that \ $\kappa\in\eta\in \OR^{\nouse^\alpha}$
\ we let \ $\mathbb{S}_\eta=(S_\gamma,\R,\kappa,\mu\cap S_\gamma)$. \ Note that \
$\mathbb{S}_\eta\in N^\alpha$ \ and the sequence \ 
$\langle \mathbb{S}_\eta\ : \ \kappa<\eta<\OR^{\nouse^\alpha} \rangle$ \ is \
$\Sigma_1(\nouse^\alpha)$.
\ Now, for \ $y\in\R$
\begin{alignat*}{1}
y\in Q&\iff \nouse^\alpha\models \varphi(y,z)\\
&\iff (\exists x\in\R)\left[\mathbb{S}_{k(x)}\models \varphi(y,z)\right]\\
&\iff \nouse^\alpha\models(\exists x\in\R)(\forall\mathcal{S})(\forall\eta)\left[(\eta=k(x) \land
\mathcal{S}=\mathbb{S}_\eta)\rightarrow \varphi(y,z)^\mathcal{S}\right]\footnotemark
\end{alignat*}
where \ $\varphi(y,z)^\mathcal{S}$ \ is the ``relativization'' of \ $\varphi$ \ to \ $\mathcal{S}$. 
\ Therefore, \ $Q\in\exists^{\R}\Pi_1(\nouse^\alpha,\R)$ \ and this completes the proof of the
Claim. 
\footnotetext{If  $k$  were not
total, then this condition would hold, vacuously, for all reals  $y$.}
\end{proof}

We conclude from the Claim that 
\[\Sigma_{n}^*(\mathfrak{B}_1(\nouse^\alpha))\subseteq\Sigma_{n}^*(\Sigma_1(\nouse^\alpha,\R))\]
and thus, 
\[\boldface{\Sigma}{n+1}(\nouse^\alpha) \subseteq \Sigma_{n}^*(\Sigma_1(\nouse^\alpha,\R)).\]
This completes the proof of the Lemma.
\end{proof}
\begin{theorem}\label{yesscale1}
Suppose that \ $\alpha$ \ properly begins a \ $\Sigma_1(\nouse)$--gap, \ $\alpha$ \ is
uncollectible, and 
\ $\nouse^{\alpha+1}\models\AD$. \ Then the pointclasses
\[
\boldface{\Sigma}{2n+1}(\nouse^\alpha) \text{ \ and \ }
\boldface{\Pi}{2n+2}(\nouse^\alpha)\]
have the scale property, for all \ $n\ge 0$.
\end{theorem}
\begin{proof} The theorem follows from Lemma~\ref{pntclsequal}, Corollary~\ref{relscales} and the
Second Periodicity Theorem of Moschovakis \cite[Theorems 6C.2 and  6C.3]{Mosch} (see Remark
\ref{spthm:rmk}).
\end{proof}

When \ $\alpha$ \ properly begins a \ $\Sigma_1(\nouse)$--gap and is collectible, then
Martin's arguments in \cite{martin} give the following analogues of Theorem 2.7 and Corollary 2.8 of
\cite{Steel}.

\begin{theorem}\label{nomore}
Suppose that \ $\alpha$ \ properly begins a \ $\Sigma_1(\nouse)$--gap,  \ $\alpha$ \ is
collectible, and \ $\nouse^{\alpha+1}\models\AD$. \ Then there is a \ $\Pi_1(\nouse^\alpha)$ \
subset of \ $\R\times\R$ \ with no uniformization in \ $N^{\alpha+1}$.
\end{theorem}

\begin{corollary}\label{noscale1}
Suppose that \ $\alpha$ \ properly begins a \ $\Sigma_1(\nouse)$--gap, \ $\alpha$ \ is collectible,
and  \ $\nouse^{\alpha+1}\models\AD$. \ Then the pointclasses
\[\boldface{\Sigma}{n+1}(\nouse^\alpha) \text{ \ and \ } \boldface{\Pi}{n}(\nouse^\alpha)\]
do not have the scale property for all \ $n\ge 1$.
\end{corollary}

\subsection{Scales inside a gap}\label{inthegap}
In this subsection we shall focus our attention on \ $\Sigma_1(\nouse)$--gaps \ $[\alpha,\beta]$ \ where \
$\alpha<\beta$. \ Thus, \ $\alpha$ \ is proper and collectible.
Our first theorem extends the above Theorem~\ref{nomore}. The proofs of Theorem 2.9 and Corollary 2.10 of
\cite{Steel} easily generalize to give the next two results. 
\begin{theorem}\label{nouniform1}
Suppose that \ the \ $\Sigma_1(\nouse)$--gap \ $[\alpha,\beta]$ \ is such that \ $\alpha<\beta$ \ 
\ and \ $\nouse^{\alpha+1}\models\AD$. \ Then there is a \ $\Pi_1(\nouse^\alpha)$ \ subset of
\ $\R\times\R$ \ with no uniformization in \ $\boldface{\Sigma}{1}(\nouse^{\beta})$.
\end{theorem}

\begin{corollary}\label{noscale2}
Suppose that \ the \ $\Sigma_1(\nouse)$--gap \ $[\alpha,\beta]$ \ is such that \ $\alpha<\beta$ \ 
\ and \ $\nouse^{\alpha+1}\models\AD$. \ If \ $\alpha<\gamma<\beta$, \ then the pointclasses
\[\boldface{\Sigma}{n}(\nouse^\gamma) \text{ \ and \ }
\boldface{\Pi}{n}(\nouse^\gamma)\]
do not have the scale property for all \ $n\ge 1$.
\end{corollary}

Thus, no new scales exist properly inside a \ $\Sigma_1(\nouse)$--gap \ $[\alpha,\beta]$. \ 
Suppose \ $\beta$ \ {\sl properly\/} ends this \ $\Sigma_1(\nouse)$--gap. 
Do any of the pointclasses \ $\boldface{\Sigma}{n}(\nouse^\beta)$ \ or \ $\boldface{\Pi}{n}(\nouse^\beta)$ \ have the scale
property? The results in
\cite{Part1} and \cite{Part2} will be used to answer this question in our next subsection.

\subsection{Scales at the proper ending of a gap}\label{endgap}
Recall that \ $\nouse=(N,\R,\kappa,\mu)$ \ is an
iterable real premouse.
In this subsection we shall again deal with \ $\Sigma_1(\nouse)$--gaps \ $[\alpha,\beta]$ \ where \ $\alpha<\beta$. \ 
The next theorem shows that the proper ending of a gap produces a real mouse.
\begin{theorem}\label{isamouse}  Suppose that \ $\beta$ \ properly ends a \ $\Sigma_1(\nouse)$--gap. Then
\ $\nouse^\beta$ \ is a real mouse.
\end{theorem}
\begin{proof} Recall that \ $\nouse=(N,\R,\kappa,\mu)$ \ is an iterable real premouse. Since \ $\nouse^\beta$ \
is an iterable real premouse, Theorem~{\GenDJ} of \cite{Part1} asserts that \ $\nouse^\beta$ \ is acceptable above the
reals. We now show that \ $\nouse^\beta$ \ is critical (see subsection~\ref{Prelims}). Because \ $\nouse^\beta$ \ is a proper initial
segment of \ $\nouse$, \ the statement \ ``$\nouse^\beta$ \ is critical'' is one whose truth can be locally verified in \ $\nouse$. \
Thus, for any premouse iterate \ $\nouse_\lambda$ \ with corresponding premouse embedding \ $\pi\colon\nouse\mapsigma{1} \nouse_\lambda$,
\ we have that \ $\nouse\models \text{``$\nouse^\beta$ \ is critical''}$ \ if and only if \ $\nouse_\lambda\models
\text{``$\pi(\nouse^\beta)$ \ is critical''}$. \ Consequently we can assume, without loss of generality, that \
$\kappa$ \ is a regular cardinal\footnote{This appeal to $\AC$ is removable. Let
$M=L(\nouse)$. One can prove the theorem in a $\ZFC$--generic extension $M[G]$  and thus, by absoluteness, the result holds in $L(\nouse)$ and hence in $V$.\label{removenote}} greater than \ $\Theta$, \ and that for all \ $X\subseteq\kappa$ \ in \ $\nouse$, \ we have \ $X\in\mu$ \ if and only if \ $X$ \ contains a closed unbounded subset of \ $\kappa$. 

We will now prove that \ $\nouse^\beta$ \ is critical. Suppose, for a contradiction, that \
$\nouse^\beta$ \ is not critical; that is, suppose that \  $\rho_{\nouse^\beta}^n>\kappa$ \  for all \ $n\in\omega$. \ For
short, let \ $\mouse=\nouse^\beta$. \ Since \ $\mouse\in\nouse$ \ it follows for each \ $n$ \ that \ $\mouse^n\in\nouse$ \
and thus, for each \ $n$ \ the truth of \ ``$\rho_{\nouse^\beta}^n>\kappa$'' can be verified in \ $\nouse$. \ Hence, for
every \ $n$, \ we have that
\ $\nouse\models \text{``$\rho_{\nouse^\beta}^n>\kappa$''}$ \ and the statement
\ ``$\rho_{\nouse^\beta}^n>\kappa$'' is equivalent to a local \ $\Sigma_1$ \ condition over \ $\nouse$. 

Because \ $\beta$ \ properly ends the \ $\Sigma_1(\nouse)$--gap,  there is a real \ $x$ \  and a \ $\Sigma_1$ \ formula \ $\varphi(v)$ \ in the
language \ $\Lng=\{\,\en,\underline{\R},\underline{\kappa},\mu\,\}$ \  such that \
$\nouse^{\beta+1}\models\varphi(x)$ \ and \ $\nouse^{\beta}\not\models\varphi(x)$. \ Thus, there is a \ $k\in\omega$ \ such
that \ $(S_{\beta+k}^{\nouse}(\R),\R,\kappa,\mu)\models\varphi(x)$ \ (see subsection~{\abovereals} of
\cite{Part1}). \ By Proposition~{\theendprop} of \cite{Part1}, there is a formula \ $\chi(v)$ \ in the language\footnote{The
language  $\Lng^\mu$  has a quantifier symbol $\mu$ (see subsection~\ref{Prelims}).\label{langnote}}
\ $\Lng^\mu$ \ such that \ $(\nouse^\beta,\mu^{\beta+1})\models\chi(x)$ \ and \
$(\nouse^\gamma,\mu^{\gamma+1})\not\models\chi(x)$  \ for all \ $\gamma<\beta$. \  Because \ $\rho_{\nouse^\beta}^n>\kappa$
\  for all \ $n\in\omega$, \ we conclude that  there is a \ $\Sigma_\omega$ \ formula\footnote{The formula $\psi(v)$ is
constructed by induction on the complexity of $\chi(v)$. See \cite[Cor. 1.33]{Cfsrm} and \cite[Cor. 2.13]{Cfsrm}.}
\ $\psi(v)$ \ in the language \ $\Lng$ \ such that  \ $\nouse^\beta\models\psi(x)$ \ and for all \ $\gamma\le\beta$, \ if \
$\nouse^\gamma\models\psi(x)$ \ then \ $(\nouse^\gamma,\mu^{\gamma+1})\models\chi(x)$.

Since \ $\rho_{\nouse^\beta}^n>\kappa$ \  for all \ $n\in\omega$, \ the proof of Theorem~{\intismi}
of \cite{Part1} applies and shows that \ $\nouse^\beta$ \ is $n$-iterable for each \ $n$. \  Now, let \ $m\in\omega$ \ be so
that \ $\chi(v)$ \ is \ $\Sigma_m$. \ Let \ $n\ge 1$ \ be fixed for a moment. Let \ $\H=\Hull_1^{\mouse^n}(\R)$.
\ Hence, \ $\H\prec_1\mouse^n$. \ Let \ $\C^n$ \ be the transitive collapse of \ $\H$. \ By Lemma~{\EOEL} of
\cite{Part1} there is a decoding \ $\C$ \ of \ $\C^n$ \ and a map \ $\sigma\colon \C \mapsigma{{n+1}} \mouse$. \ 
Hence, \ $\C$ \ is an iterable
premouse and is acceptable above the reals. It thus follows, from the definition of \ $\H$, \ that \ $\C$ \ is
critical and so, \ $n(\C)$ \ is defined. Now, assume that \ $n$ \ is large enough to ensure that \ $\n=n(\C)>m$.
\ This assumption implies that \ $\C\models\chi(x)$  \ and \
$\C\models\lnot\varphi(x)$. \  Theorem
{\criterion} of \cite{Part1} implies that \ $\C$ \ is \ $\n$--iterable. Therefore, \ $\C$ \ is a real mouse.

Let \ $\C_\kappa$ \ be the \ $\kappa^{\,\ul{\text{th}}}$ \ mouse iterate of \
$\C$. \ Since \ $\kappa>\Theta$, \ our assumption on \ $\kappa$ \ implies that \ $\C_\kappa$ \  and \
$\nouse$ \ are comparable (see \cite[Definition
2.23]{Crcm}).  

\begin{claim} $\C_\kappa=\nouse^\beta$.
\end{claim}
\begin{proof} Since \ $\C\models\lnot\varphi(x)$, \ we see that \ $\C_\kappa\models\lnot\varphi(x)$. \ If \
$\nouse^\beta$
\ were a proper initial segment of \
$\C_\kappa$, \ then  since \
$\nouse^\beta\models\varphi(x)$, \ we would conclude that \ $\C_\kappa\models\varphi(x)$, \ as \
$\varphi(v)$ \ is \ $\Sigma_1$. \ But this is not possible, because \ $\C_\kappa\models\lnot\varphi(x)$. \ Thus, we must
have that \
$\C_\kappa=\nouse^\gamma$ \ for some \ $\gamma\le\beta$. \ Since \ $\C\models\chi(x)$ \ and \ $\chi(v)$ \
is \
$\Sigma_m$, \ it follows that \ $\C_\kappa \models\chi(x)$ \ because \ $n(\C)>m$. \ Consequently, \
$\nouse^\gamma\models\chi(x)$ \ and, as noted above,  we have that \
$(\nouse^\gamma,\mu^{\gamma+1})\models\psi(x)$. \ Hence, by Proposition~{\theendprop} of \cite{Part1} we have  \ $\nouse^{\gamma+1}\models\varphi(x)$. \ Therefore, \
$\nouse^\gamma=\nouse^\beta$ \ and this completes the proof of the claim.
\end{proof}

Since a mouse iterate of a mouse is again a mouse, we conclude from the Claim that \ $\nouse^\beta$ \  must be
critical. This contradiction forces us to conclude that if \
$\beta$ \ properly ends a \
$\Sigma_1(\nouse)$--gap, then \ $\nouse^\beta$ \ is critical. \ Theorem~{\intismi} of \cite{Part1} now implies that \
$\nouse^\beta$ \ is a real mouse.
\end{proof}

\begin{theorem}\label{newset}  Suppose that \ $\beta$ \ properly ends a \ 
$\Sigma_1(\nouse)$--gap,
\ where \ $\nouse=(N,\R,\kappa,\mu)$, \ and let \ $\n=n(\nouse^\beta)$. \ Then \ $(\exists
A\subseteq\R)\left(A\notin N^{\beta+1}\setminus N^\beta\right)$ \ if and only if \ $\rho^{\n+1}_{\nouse^\beta}<\kappa$.
\end{theorem}
\begin{proof} Suppose that \ $\beta$ \ properly ends a \ 
$\Sigma_1(\nouse)$--gap,
\ where \ $\nouse=(N,\R,\kappa,\mu)$. \ Theorem~\ref{isamouse} implies that \ $\nouse^\beta$ \ is a real
mouse and so, \ $\n=n(\nouse^\beta)$ \ is defined. Thus, \ $\rho^{\n+1}_{\nouse^\beta}\le\kappa$. \ 
If there exists a new set of reals in \ $\nouse^{\beta+1}$ \ which is not in \ $\nouse^\beta$, \ then Lemma
{\lemmatwo} of \cite{Part1} implies that \ $\rho^{\n+1}_{\nouse^\beta}<\kappa$. 

Now assume that \ $\rho^{\n+1}_{\nouse^\beta}<\kappa$. \ We prove that there is a set of reals \ $A$ \ in \ $\nouse^{\beta+1}$ \ which is not in \
$\nouse^\beta$. \ As  in the proof of the previous theorem, we shall assume\footnote{See footnote~\ref{removenote}.} that \
$\kappa$ \ 
is a regular cardinal greater than
\ $\Theta$ \ and that for all \ $X\subseteq\kappa$ \ in \ $\nouse$ \ we have \ $X\in\mu$ \ if and only if \ $X$ \ contains a closed unbounded subset of \ $\kappa$. 
\ Let \ $\core=\core(\nouse^\beta)$ \ be the core of the
real mouse \ $\nouse^\beta$ \ (see subsection~{\coremice} of \cite{Part1}), let \ $\kappa_0=\kappa^\core$, \
$\overline{\core}=\core^\n$ \ and let \ $\overline{\nouse^\beta}=(\nouse^\beta)^{\n}$. \ Recall that \
$\overline{\core}$ \ and \ $\overline{\nouse^\beta}$ \ are structures for the language \ $\Lng_\n$ \ (see subsection
\ref{Prelims}). The real mouse \ $\nouse^\beta$ \ is a mouse iterate of its core \ $\core$, \ and \ $\overline{\nouse^\beta}$ \
is a premouse iterate of \ $\overcore$. \ In particular, Lemma~\ref{deftwotwo}  implies that the mouse iterate \
$\core_\kappa=\nouse^\beta$ \ and the premouse iterate \ ${\overcore}_\kappa=\overline{\nouse^\beta}$. \ Let \
$\overline{\pi}_{0\kappa}\colon\overcore\mapsigma{1}{\overline{\nouse^\beta}}$ \ be the premouse embedding of \ $\overcore$ \ into
its \ $\kappa^{\,\ul{\text{th}}}$ \ premouse iterate \ $\overline{\nouse^\beta}$.

\begin{claim} $\rho_{\overline{\core}}^{k+1}=1$ and \ $\rho_{\overline{\core}}^{k}>1$ \ for some \ $k\ge 0$.
\end{claim}
Before we prove this claim we show how it can be used to complete the proof of the theorem. Let \ $k\in\omega$
\ be as in the statement of the Claim. Corollary~{\gammaproj} of \cite{Part1} (also see \cite[Definition
{\defproj}]{Part1}) implies that there exists  a \ $\boldface{\Sigma}{k+1}(\overline{\core})$ \ set of reals \ $A$ \ 
not in \ $\overline{\core}$. \ Because the real mouse \ $\nouse^\beta$ \ is a mouse iterate of its core \ $\core$ \ (see Theorem
2.33 of \cite{Cfsrm}), \ Lemma 2.19 of \cite{Cfsrm} implies that \ $A$ \ is a
\ $\boldface{\Sigma}{\n+k+1}(\nouse^\beta)$  \ set of reals that is not in \ $\nouse^\beta$. \ Hence, the desired conclusion
follows. 
\begin{proof}[Proof of Claim] Assume, for a contradiction, that \ $\rho_{\overline{\core}}^{k}>1$ \ for all \ $k\ge 1$. \ 
Since \ $\beta$ \ properly ends the \ $\Sigma_1(\nouse)$--gap,  there is a \ $\Sigma_1$ \ formula \ $\varphi(v)$ \ in the language \
$\Lng=\{\,\en,\underline{\R},\underline{\kappa},\mu\,\}$ \ and a real \ $x$ \ such that \ $\nouse^{\beta+1}\models\varphi(x)$ \
and \ $\nouse^{\beta}\not\models\varphi(x)$. \ Thus, there is an \ $i\in\omega$ \ such that \
$(S_{\beta+i}^{\nouse}(\R),\R,\kappa,\mu)\models\varphi(x)$ \ (see subsection~{\abovereals} of
\cite{Part1}). \ By Proposition~{\theendprop} of \cite{Part1}, there is a formula \ $\chi(v)$ \ in the language\footnote{See
footnote~\ref{langnote}.} \ $\Lng^\mu$ \ such that \ $(\nouse^\beta,\mu^{\beta+1})\models\chi(x)$ \ and \
$(\nouse^\gamma,\mu^{\gamma+1})\not\models\chi(x)$  \ for all \ $\gamma<\beta$. \  Because \
$\rho_{\nouse^\beta}^{\n+1}<\kappa$, \ Lemma~\ref{deftwotwo} asserts that \ $\mu^{\beta+1}$ \ is
$\kappa_0$--predictable (see Definition~\ref{preddefagain}).  Therefore, as in the proof of Lemma~{\predlemma} of
\cite{Part1}, there is a \ $\Sigma_\omega$ \ formula \ $\psi(u,v)$ \ in the language \ $\Lng$ \ such that  
\begin{equation}(\nouse^\beta,\mu^{\beta+1})\models\chi(x)\iff \nouse^\beta\models\psi(x,\kappa_0)\label{wacky1}\end{equation}
where \ $\psi$ \ depends only on \ $\chi$ \ and \ $\n$; \ that is, the equivalence (\ref{wacky1}) holds for any such real mouse \
$\nouse^\beta$ \ where \ $\kappa_0$ \ is the ``measurable cardinal'' of its core. 

It follows from Theorem~{\MET} of
\cite{Part1} (also see \cite[Theorem 2.11]{Cfsrm}) \ that there is a canonical \ $\Sigma_k$ \ formula \ $\overline{\psi}(u,v)$ \
in the language \ $\Lng_\n$ \ such that 
\begin{equation}\nouse^\beta\models\psi(x,\kappa_0)\iff
\overline{\nouse^\beta}\models\overline{\psi}(x,\kappa_0).\label{wacky2}\end{equation}
Let \ $f\in \core$ \ be the identity function \
$f\colon\kappa_0\to\kappa_0$. \ Because \ $\overline{\pi}_{0\kappa}(f)(\kappa_0)=\kappa_0$, \ 
Lemma 2.19 of \cite{Crcm} implies that there is a \ $\Sigma_k$ \ formula \ $\overline{\psi}^*(u,v)$ \ in the language \
$\Lng_\n$, \ depending only on \ $\overline{\psi}$, \ so that  \ 
\begin{equation}\overline{\nouse^\beta}\models\overline{\psi}(x,\kappa_0)\iff
\overcore\models\overline{\psi}^*(x,f).\label{wacky3}\end{equation}
Lemma~{\soundnplusone} and its Corollary~{\skolem} of \cite{Part1} assert that
\ $\overline{\core}$ \ is \ $(k+1)$--sound and \ $\overline{\core}$ \ has a \ $\Sigma_{k+1}$ \ Skolem function. Let \ $a$ \ be
an element of \ $\overcore$ \ so that \ $\overcore$ \ has a \ $\Sigma_{k+1}$ \ Skolem function which is \ $\Sigma_{k+1}(\{a\})$ \
(see \cite[Definition~{\skolemdefn}]{Part1}). Let \ $\mathcal{H} = \Hull_{k+1}^{\overcore}(\R\cup \{a, f\})$, \ and let  
\ $\overkore$ \  be the transitive collapse of \ $\mathcal{H}$. \ Let \ $\overline{\sigma}\colon\overkore\mapsigma{k+1}\overcore$
\ be the inverse of the collapse map. As in the proof of Lemma 2.29 of \cite{Cfsrm}, there is a real mouse \ 
$\kore$ \ such that \ $n(\kore)=n(\core)=\n$ \ and \ $\rho^{\n+1}_{\kore}\le\rho^{\n+1}_{\core}<\kappa$. \ Note that \ $\kore$ \
is a core mouse \ where \ $\overline{\sigma}(\kappa_0)=\kappa^\kore$ \ and   \
$\overline{\sigma}(f)\colon\overline{\sigma}(\kappa_0)\to\overline{\sigma}(\kappa_0)$ \ is the identity function in \
$\overkore$. \ In addition,
\ $\overkore\models\overline{\psi}^*(x,\overline{\sigma}(f))$. \ Since \ $\kappa$ \ is a regular cardinal greater than
\ $\Theta$, \ it follows that the mouse iterates \ $\core_\kappa$ \ and \ $\kore_\kappa$ \ are comparable.  Due to the fact
that
\ $\overline{\sigma}\colon\overkore\mapsigma{k+1}\overcore$, \ we conclude that \ $\kore_\kappa$ \ must be an initial segment
of \ $\core_\kappa$. \ Since \ $\core_\kappa=\nouse^\beta$, \ there must be an ordinal \ $\gamma\le\beta$ \ such that \
$\kore_\kappa=\nouse^\gamma$. \ Thus, \ $\nouse^\gamma$ \ is a real mouse with core \ $\kore$. \ Recall that \
$\overkore\models\overline{\psi}^*(x,\overline{\sigma}(f))$. \ Because the biconditionals (\ref{wacky1})--(\ref{wacky2}) are
sufficiently uniform, we conclude that \ $(\nouse^\gamma,\mu^{\gamma+1})\models\chi(x)$ \ and hence,
\ $\nouse^{\gamma+1}\models\varphi(x)$. \ Therefore, \
$\nouse^{\beta+1}=\nouse^{\gamma+1}$ \ and thus, \ $\nouse^{\beta}=\nouse^{\gamma}$. \ Consequently, \ $\core=\kore$ \ and \
$\overcore=\overkore$. \ It now follows from the construction of \ $\overkore$ \ that there is a  \
$\boldface{\Sigma}{k+1}(\overcore)$ \ set of reals not in \ $\overcore$. \ Therefore, \ $\rho_{\overline{\core}}^{k+1}=1$ \ and
this contradiction ends our proof of the Claim.
\end{proof}

The proof of the theorem is now complete.\end{proof}

\begin{remark} Suppose that \ $\beta$ \ properly ends a \ $\Sigma_1(\nouse)$--gap and let \ $\n=n(\nouse^\beta)$.
\ Lemma~{\deftwo}(5) of \cite{Part1} and the above Theorem~\ref{newset} imply that \ $m(\nouse^{\beta})$ \ is
defined if and only if  \ $\rho^{\n+1}_{\nouse^\beta}<\kappa$.
\end{remark}

Suppose that \ $\beta$ \ properly ends a \ $\Sigma_1(\nouse)$--gap.
Theorem~\ref{newset} implies that if \ $\pow(\R)\cap N^{\beta+1}\setminus N^\beta=\emptyset$, \
then \ $\nouse^{\beta+1}$ \ is ``$\R$--sharplike''\footnote{One can generalize Dodd's
theorems (see \cite[Chapter 15]{Dodd}) concerning ``sharplike mice'' to encompass ``$\R$--sharplike real
mice''. \label{sharpnote}}. It follows that $\R$--sharplike real mice relate to an inner model of \ $V\eq\Kr$ \ in the
same way that \ $\R^\shrp$ \ relates to \ $\Lr$. \ This, together with Theorem 5.17 of \cite{Crcm}, allows one to prove the
following theorem.

\begin{theorem}
Let \ $\nouse=(N,\R,\kappa,\mu)$ \ be an iterable real premouse. Suppose that
\ $\beta$ \ properly ends a \ $\Sigma_1(\nouse)$--gap and  let \
$\n=n(\nouse^\beta)$. \ Assuming \ $\nouse^{\beta}\models\AD$, \ if \
$\rho_{\nouse^\beta}^{\n+1}=\kappa$  \ then the pointclasses
\ $\boldface{\Sigma}{n}(\nouse^\beta)$  \ and \ $\boldface{\Pi}{n}(\nouse^\beta)$ \ 
do not have the scale property for all \ $n\ge 1$.
\end{theorem}
\begin{proof}[Sketch of Proof] For \ $\gamma\le \widehat{\OR}^{\,\nouse}$ \ let \ $H^{\nouse^\gamma}_{\kappa}=\{a\in
N^\gamma : 
\abs{T_c(a)}_{\nouse^\gamma}<\kappa\}$ \ (see
\cite[Definition
{\trans}]{Part1}). Since \ $\rho_{\nouse^\beta}^{\n+1}=\kappa$, \ Lemma~{\lemmatwo} of \cite{Part1}) implies that \
$H^{\nouse^\beta}_{\kappa}=H^{\nouse^{\beta+1}}_{\kappa}$. \  One can prove (see
\cite[Chapter 15]{Dodd}) that the sets of reals in \ $H^{\nouse^\beta}_{\kappa}$ \ are exactly those in an inner model of \
$V\eq\Kr$. \ Thus, Theorem 5.17 of \cite{Crcm} implies that there is a set of reals in \ $H^{\nouse^\beta}_{\kappa}$ \ which
has no scale in \ $H^{\nouse^\beta}_{\kappa}$. \ Because \
$H^{\nouse^\beta}_{\kappa}=H^{\nouse^{\beta+1}}_{\kappa}$,
\ it follows  that \
$\boldface{\Sigma}{n}(\nouse^\beta)=H^{\nouse^\beta}_{\kappa}$ \ as pointclasses, for each \ $n$. \ Therefore, \
$\boldface{\Sigma}{n}(\nouse^\beta)$ \ and \ $\boldface{\Pi}{n}(\nouse^\beta)$
\ do not have the scale property.
\end{proof}

\begin{lemma}\label{subsetrel2}  Let \ $\nouse=(N,\R,\kappa,\mu)$ \ be an iterable real premouse and
suppose that \ $\beta$ \ properly ends a \ $\Sigma_1(\nouse)$--gap and  let \ $\n=n(\nouse^\beta)$. \ Suppose that
\ $\rho_{\nouse^\beta}^{\n+1}<\kappa$ \ and let \ $\m=m(\nouse^{\beta})$. \ Then 
\[\boldface{\Sigma}{n+\m}(\nouse^\beta)\subseteq
\Sigma_{n}^*(\mathfrak{B}(\boldface{\Sigma}{\m}(\nouse^\beta)))\]
as pointclasses, for all \ $n\ge 0$.
\end{lemma}
\begin{proof} Because \ $\rho_{\nouse^\beta}^{\n+1}<\kappa$, \ Theorem~\ref{newset} implies that \
$\m=m(\nouse^\beta)$ \ is defined. Let \ $\core=\core(\nouse^\beta)$. \ Because \ $\nouse^\beta$ \ is a mouse iterate of \
$\core$, \ we see that \ $m(\core)=\m$ \ by Lemma 2.19 of \cite{Cfsrm}. Note that \
$\rho_{\nouse^\beta}^{0}>1$  \ by definition (see \cite[Definition 1.18]{Cfsrm}). We will now show that \
$\rho_{\nouse^\beta}^{i}>1$ \ for all \ $0\le i<\m$. \ Suppose that for some \ $i<\m$ \ we have that \
$\rho_{\nouse^\beta}^{i}=1$. \ Assume that \ $i\ge 1$ \ is the smallest such natural number. Thus, \
$\rho_{\nouse^\beta}^{i-1}>1$. \ Corollary 2.38 of \cite{Cfsrm} then implies that there is a \
$\boldface{\Sigma}{i}(\nouse^\beta)$ \ set of reals not in \ $\nouse^\beta$, \ which is impossible because \ $i<\m$.
\ We conclude \ $\rho_{\nouse^\beta}^{\m-1}>1$ \ and thus, \ Corollary 2.38 of
\cite{Cfsrm} \ also implies that \ $\rho_{\core}^{\m-1}>1$. \ By Lemma 2.34 of \cite{Cfsrm} and Corollary 1.32  
of \cite{Cfsrm},  there is a partial \ $\boldface{\Sigma}{\m}(\core)$ \ map \ $g\colon\R\maps{onto}C$, \ where \ $C$ \ is
the domain of \ $\core$. \ For \ $n\ge 0$ \ it follows, as in Proposition~\ref{propp},\footnote{See Examples~\ref{exmone} and
\ref{exmtwo}.} that \ $\boldface{\Sigma}{n+\m}(\core)\subseteq\Sigma_{n}^*(\mathfrak{B}(\boldface{\Sigma}{\m}(\core)))$ \ as
pointclasses. We know that \ $\nouse^\beta$ \  is a mouse iterate of \ $\core$ \ and thus,
\ $\boldface{\Sigma}{k+\m}(\core)=\boldface{\Sigma}{k+\m}(\nouse^\beta)$ \ as pointclasses for \ $k= 0$ \ and \ $k=n$, \ by
Lemma 2.19 of \cite{Cfsrm}. Therefore, \
$\boldface{\Sigma}{n+\m}(\nouse^\beta)\subseteq \Sigma_{n}^*(\mathfrak{B}(\boldface{\Sigma}{\m}(\nouse^\beta)))$.
\end{proof}
\begin{definition} Let \ $\nouse$ \ be an iterable real premouse and let \
$[\alpha,\beta]$ \ be a \ $\Sigma_1(\nouse)$--gap. If \
$\alpha<\beta<\widehat{\OR}^{\,\nouse}$ \ and \ $m(\nouse^{\beta})$ \ is defined,  
then we shall say that 
\begin{itemize}
\item $\beta$ \ is {\it weak\/} if and only if \ $\nouse^\beta$ \ is a weak mouse,
\item $\beta$ \ is {\it strong\/} if and only if \ $\nouse^\beta$ \ is a strong mouse.
\end{itemize}
\end{definition}

Our next theorem is an observation that can be used  to prove within \ $\Kr$ \ that the axiom of determinacy is equivalent to the existence of arbitrarily large cardinals \ $\kappa<\Theta$ \ with the strong partition property  (see \cite[page 432]{Kana}).
\begin{theorem}\label{scales2}  Suppose that \ $[\alpha,\beta]$ \ is a \
$\Sigma_1(\nouse)$--gap where  \ $\beta$ \ is weak. Let \ $\m=m(\nouse^\beta)$. \ If \
$\nouse^{\alpha}\models\AD$, \ then the pointclass \
$\boldface{\Sigma}{\m}(\nouse^\beta)$ \ has the scale property.
\end{theorem}
\begin{proof} Suppose that \ $\nouse^{\alpha}\models\AD$. \ If \ $\nouse^\beta\models\AD$, \ then Theorem~\ref{newthmrpt} 
asserts that \ $\boldface{\Sigma}{\m}(\nouse^\beta)$
\ has the scale property. We shall prove that \ $\nouse^\beta\models\AD$. \ Assume, for a contradiction, that  \
$\nouse^\beta\not\models\AD$. \ Thus, there is a non-determined game in \ $\nouse^\beta$ \ and this can be asserted as a  \ $\Sigma_1$ \
statement true in \ $\nouse^\beta$. \ Since \ $\nouse^{\alpha}\prec_1^{\R}\nouse^{\beta}$, \ there is a non-determined game in \
$\nouse^\alpha$. \ This contradiction completes the proof of the theorem.
\end{proof}

\begin{lemma}\label{pntclsequal2}  Suppose that \ $[\alpha,\beta]$ \ is a \ $\Sigma_1(\nouse)$--gap
where \ $\beta$  \ is weak. Let \ $\m=m(\nouse^\beta)$. \ Then 
\begin{gather}
\boldface{\Sigma}{n+\m}(\nouse^\beta) = \Sigma_{n}^*(\boldface{\Sigma}{\m}(\nouse^\beta))\tag{a}\\
\boldface{\Pi}{n+\m}(\nouse^\beta) = \Pi_{n}^*(\boldface{\Sigma}{\m}(\nouse^\beta))\tag{b}
\end{gather}
as pointclasses,  for all \ $n\ge 0$.
\end{lemma}
\begin{proof} For \ $n=0$, \ the conclusion holds by  Definition~\ref{realquantifier}. So we assume \ $n\ge 1$. \ It is
sufficient to show that \ $\boldface{\Sigma}{n+\m}(\nouse^\beta) =
\Sigma_{n}^*(\boldface{\Sigma}{\m}(\nouse^\beta))$. \ Clearly, \
$\Sigma_{n}^*(\boldface{\Sigma}{\m}(\nouse^\beta))\subseteq\boldface{\Sigma}{n+\m}(\nouse^\beta)$. \ We show that \
$\boldface{\Sigma}{n+\m}(\nouse^\beta)\subseteq\Sigma_{n}^*(\boldface{\Sigma}{\m}(\nouse^\beta))$. \ By Lemma~\ref{subsetrel2}
we have that \[\boldface{\Sigma}{n+\m}(\nouse^\beta)\subseteq
\Sigma_{n}^*(\mathfrak{B}(\boldface{\Sigma}{\m}(\nouse^\beta)))\]
and thus it is enough to show, as in the proof of Lemma~\ref{pntclsequal}, that \ $\boldface{\Sigma}{\m}(\nouse^\beta)\subseteq
\exists^\R\boldface{\Pi}{\m}(\nouse^\beta)$. \ To do this, let \ $P\subseteq\R$  \ be in \
$\boldface{\Sigma}{\m}(\nouse^\beta)$. \ Because \ $\beta$ \ is weak, Theorem~\ref{essential} states
that there is a total
\ $\boldface{\Sigma}{\m}(\nouse^\beta)$ \  map \ $k\colon\omega\to N^\beta$ \ such that \ 
$P=\bigcup\limits_{i\in\omega}k(i)$ \ and thus
\[P(x) \iff (\exists i\in\omega)(\forall a)(a=k(i)\rightarrow x\in a)\]
for all \ $x\in\R$. \ Since the graph of \ $k$ \ is  \ $\boldface{\Sigma}{\m}(\nouse^\beta)$, \ we conclude that \ $P$ \ is in
\ $\exists^\omega\boldface{\Pi}{\m}(\nouse^\beta)$ \ and so, \ $P$ \ is in  \ $\exists^\R\boldface{\Pi}{\m}(\nouse^\beta)$.
\end{proof}
\begin{theorem}\label{yesscale2}Suppose that \ $[\alpha,\beta]$ \ is a \ $\Sigma_1(\nouse)$--gap and \ $\beta$ \ properly 
ends this gap. If \ $\beta$ \ is weak, \ $\m=m(\nouse^\beta)$ \ and \ $\nouse^{\beta+1}\models\AD$, \ then the
pointclasses
\[
\boldface{\Sigma}{\m+2k}(\nouse^\beta) \text{ \ and \ }
\boldface{\Pi}{\m+(2k+1)}(\nouse^\beta)\]
have the scale property, for all \ $k\ge 0$.
\end{theorem}
\begin{proof} This follows directly from Lemma~\ref{pntclsequal2}, Theorem~\ref{scales2} and the
Second Periodicity Theorem of Moschovakis \cite[Theorems 6C.2 and  6C.3]{Mosch} (see Remark
\ref{spthm:rmk}).
\end{proof}

When \ $\beta$ \ properly begins a \ $\Sigma_1(\nouse)$--gap and is strong, then
Martin's arguments in \cite{martin} also give the following analogues of Theorem 3.3 and Corollary 3.4, respectively, in
\cite{Steel}.

\begin{theorem}\label{nouniform2}
Suppose that \ $[\alpha,\beta]$ \ is a \ $\Sigma_1(\nouse)$--gap and \ $\beta$ \ properly ends this gap. If \ $\beta$
is strong  and \ $\nouse^{\beta+1}\models\AD$,  \ then there is a \
$\Pi_1(\nouse^\alpha)$ \ subset of \ $\R\times\R$ \ with no uniformization in \ $N^{\beta+1}$.
\end{theorem}

\begin{corollary}\label{noscale3}
Suppose that \ $[\alpha,\beta]$ \ is a \ $\Sigma_1(\nouse)$--gap and \ $\beta$ \ properly ends this gap. If \ $\beta$
is strong  and \ $\nouse^{\beta+1}\models\AD$,  \ then the pointclasses
\[\boldface{\Sigma}{n+1}(\nouse^\beta) \text{ \ and \ }
\boldface{\Pi}{n}(\nouse^\beta)\]
do not have the scale property for all \ $n\ge 1$.
\end{corollary}

\subsection{The scale table}
In this subsection we shall assume that the iterable real premouse \ $\nouse=(N,\R,\kappa,\mu)$
\ is such that \ $\nouse\models\AD$ \ and thus, \ $\Sigma_1(\nouse)$ \
has the scale property. In the previous subsections we identified precisely those internal levels
of the Levy hierarchy for \ $\nouse$ \ which also possess the scale property. Table
\ref{scaletable} presents a summary of this development where \ $[\alpha,\beta]$ \ is a \
$\Sigma_1(\nouse)$--gap, \ $n\ge 1$, \ $\O=\widehat{\OR}^{\,\nouse}$, \ $\n=n(\nouse^\beta)$ \ and \ $\m=m(\nouse^\beta)$ \
whenever \ $m(\nouse^\beta)$ \ is defined. In Table~\ref{scaletable}, items 3--6 focus on the proper beginning of a
$\Sigma_1(\nouse)$--gap, items 7--8 address the interior of such a gap and items 9--14 concentrate on the proper ending of a
gap. When \
$\beta$ \ is not proper, then \ $\nouse = \nouse^\beta$ \ and Table~\ref{scaletable} does not address the question of whether or
not the external pointclasses \ $\boldface{\Sigma}{n}(\nouse)$ \ or \
$\boldface{\Pi}{n}(\nouse)$ \ have the scale property for arbitrary \ $n$. \ In Section~\ref{questionQ}, we shall pursue
this issue. 

\begin{table}
\renewcommand{\arraystretch}{1.45}
\[
\begin{array}{|c|c|c|c|}\hline
\multicolumn{1}{|c}{}&\multicolumn{1}{|c|}{\textup{\sc Pointclass}}
&\multicolumn{1}{c|}{\textup{\sc Gap Property}}&\multicolumn{1}{c|}{\textup{\sc Scale Property}}\\
\hline\hline
\newresult&\boldface{\Sigma}{1}(\nouse^\alpha)&{\alpha\le\O}&\textup{ Yes }\\
\hline
\newresult&\boldface{\Pi}{1}(\nouse^\alpha)&{\alpha\le\O}&\textup{ No }\\
\hline
\newresult&\boldface{\Sigma}{n}(\nouse^\alpha)&\alpha\textup{ is uncollectible \& $\alpha<\O$}&\textup{Yes iff
$n$ is odd}\\
\hline
\newresult&\boldface{\Pi}{n}(\nouse^\alpha)&\alpha\textup{ is uncollectible \& $\alpha<\O$}&\textup{Yes iff
$n$ is even}\\
\hline
\newresult&\boldface{\Sigma}{n+1}(\nouse^\alpha)&\alpha\textup{ is collectible \& $\alpha<\O$}&\textup{No}\\
\hline
\newresult&\boldface{\Pi}{n}(\nouse^\alpha)&\alpha\textup{ is collectible \& $\alpha<\O$} &\textup{No}\\
\hline
\newresult&\boldface{\Sigma}{n}(\nouse^\gamma)&\alpha<\gamma<\beta&\textup{ No }\\
\hline
\newresult&\boldface{\Pi}{n}(\nouse^\gamma)&\alpha<\gamma<\beta&\textup{ No }\\
\hline
\newresult&\boldface{\Sigma}{n}(\nouse^\beta)&\beta\textup{\, is weak \& $\beta<\O$} &\textup{Yes iff
$(n\!-\!\m)\!\ge\!0$ is even}\\
\hline
\newresult&\boldface{\Pi}{n}(\nouse^\beta)&\beta\textup{\, is weak \& $\beta<\O$} &\textup{Yes iff
$(n\!-\!\m)\!\ge\!1$ is odd}\\
\hline
\newresult&\boldface{\Sigma}{n}(\nouse^\beta)&\beta\textup{\, is strong \& $\beta<\O$} &\textup{ No }\\
\hline
\newresult&\boldface{\Pi}{n}(\nouse^\beta)&\beta\textup{\, is strong \& $\beta<\O$}&\textup{ No }\\
\hline
\newresult&\boldface{\Sigma}{n}(\nouse^\beta)&\rho_{\nouse^\beta}^{\n}=\kappa\textup{\, \&
$\beta<\O$}&\textup{ No }\\
\hline
\newresult&\boldface{\Pi}{n}(\nouse^\beta)&\rho_{\nouse^\beta}^{\n}=\kappa\textup{\, \&
$\beta<\O$}&\textup{ No }\\
\hline
\end{array}
\]
\caption{Scale analysis of a $\Sigma_1(\nouse)$--gap $[\alpha,\beta]$ where $\O=\widehat{\OR}^{\,\nouse}$, 
$\n=n(\nouse^\beta)$ and $\m=m(\nouse^\beta)$.} 
\label{scaletable}
\end{table}

This completes the results of Section~\ref{minimal}. It turns out that these results will allow us to answer
all of our questions concerning the complexity of scales in the inner model \ $\Kr$ \ (see Section
\ref{minimalKr}).
\section{Premouse iteration preserves {\mathversion{bold}$\Sigma_1$}--gaps}\label{extensions} 

In this section we shall show the premouse iteration preserves $\Sigma_1$--gaps and preserves internal pointclasses. Then we
will show that the  iteration of a real 1--mouse also preserves its external pointclasses. Throughout this section \ $\nouse$ \
will be an iterable premouse with premouse iteration
$$\prenousesystem.$$
Clearly, \
$\widehat{\OR}^{\,\nouse}\le{\OR}^{\,\nouse}$. \ If
\ $\widehat{\OR}^{\,\nouse}={\OR}^{\,\nouse}=\delta$, \ it will be convenient to extend the domain
of each \ $\pi_\gamma$ \ to include \
$\delta$ \ by defining \
$\pi_\gamma(\delta)=\sup\left\{\pi_\gamma(\lambda) : \lambda<\delta\right\}$. 

\begin{theorem}\label{elem} Let \ $\nouse$ \ be an iterable premouse and let \ $\gamma$ \ be an ordinal.
For all \ $\alpha<\widehat{\OR}^{\,\nouse}$, \ we have that  \ $\boldface{\Sigma}{n}(\nouse^\alpha)=
\boldface{\Sigma}{n}\!\left(\nouse^{\,\pi_\gamma(\alpha)}_\gamma\right)$ \ 
as pointclasses, for all \ $n\ge1$. 
\end{theorem}
\begin{proof} Let \ $n\ge1$.  \ Because \ $\alpha<\widehat{\OR}^{\,\nouse}$ \ and \
$\pi_\gamma\colon\nouse\mapsigma{1}\nouse_\gamma$, \ we conclude that 
\ $\pi_\gamma:\nouse^\alpha\mapsigma{\omega}\nouse^{\,\pi_\gamma(\alpha)}_\gamma$. \ Consequently, we have that \ 
$\boldface{\Sigma}{n}(\nouse^\alpha)\subseteq
\boldface{\Sigma}{n}\!\left(\nouse^{\,\pi_\gamma(\alpha)}_\gamma\right)$, \ as pointclasses. For the other direction, suppose that \ $A$
\ is a set of reals in \ $\boldface{\Sigma}{n}\!\left(\nouse^{\,\pi_\gamma(\alpha)}_\gamma\right)$. \ Let \ $\varphi(u,v)$ \ be a \
$\Sigma_n$ \ formula and let \ $c$ \ be an element in \ $\nouse^{\,\pi_\gamma(\alpha)}$ \ so that 
\[x\in A \iff \nouse^{\,\pi_\gamma(\alpha)}_\gamma\models\varphi(x,c)\]
for all \ $x\in\R$. \ Let \ $N^{\pi_\gamma(\alpha)}_\gamma$ \ be the domain of \ $\nouse^{\,\pi_\gamma(\alpha)}_\gamma$. \ Since \ $A$ \
and \
$N^{\,\pi_\gamma(\alpha)}_\gamma$ \ are in
\ $\nouse_\gamma$, \ it follows that 
\[\nouse_\gamma\models (\exists y\in N^{\pi_\gamma(\alpha)}_\gamma)(\forall x\in\R)[x\in A \leftrightarrow
\nouse^{\,\pi_\gamma(\alpha)}_\gamma\models\varphi(x,y)].\]
Because \ $A$ \ is also in \ $\nouse$ \ and \ $\pi_\gamma\colon\nouse\mapsigma{1}\nouse_\gamma$, \ we see that
\[\nouse\models (\exists y\in N^{\alpha})(\forall x\in\R)[x\in A \leftrightarrow
\nouse^{\,\alpha}\models\varphi(x,y)]\]
where \ $N^{\alpha}$ \ is the domain of \ $\nouse^\alpha$. \ Therefore, \ $A$ \ is in \ $\boldface{\Sigma}{n}(\nouse^\alpha)$.
\end{proof}

\begin{theorem}\label{gap_preseve} Let \ $\nouse$ \ be an iterable premouse and let \ $[\alpha,\beta]$ \ be
a \ $\Sigma_1(\nouse)$--gap. For each \ $\gamma\in\OR$,
\ $[\pi_\gamma(\alpha),\pi_\gamma(\beta)]$ \ is a \ $\Sigma_1(\nouse_\gamma)$--gap. \ In addition, \
$\boldface{\Sigma}{1}(\nouse^\alpha)=\boldface{\Sigma}{1}\!\left(\nouse^{\,\pi_\gamma(\alpha)}_\gamma\right)$ \ as pointclasses.
\end{theorem}
\begin{proof} Let \ $[\alpha,\beta]$ \ be
a \ $\Sigma_1(\nouse)$--gap. Because \ $\pi_\gamma\colon\nouse\updownmap{cofinal}{1}\nouse_\gamma$ \ it follows easily that \
$[\pi_\gamma(\alpha),\pi_\gamma(\beta)]$ \ is a \ $\Sigma_1(\nouse_\gamma)$--gap. If \ $\alpha$ \ is not proper, then Corollary
2.14(2) of \cite{Crcm} implies that \
$\boldface{\Sigma}{1}(\nouse^\alpha)=\boldface{\Sigma}{1}\!\left(\nouse^{\,\pi_\gamma(\alpha)}_\gamma\right)$ \ as pointclasses.
If \ $\alpha$ \ is proper, then Theorem~\ref{elem} implies the desired conclusion.\end{proof}
\begin{corollary} Let \ $\nouse$ \ be an iterable premouse and let \ $\gamma$ \ be any ordinal.
\be
\item If \ $\alpha$ \ properly begins a \ $\Sigma_1(\nouse)$--gap, then \
$\pi_\gamma(\alpha)$ \ properly begins a \ $\Sigma_1(\nouse_\gamma)$--gap. 
\item If \ $\beta$ \ properly
ends a \ $\Sigma_1(\nouse)$--gap, then \ $\pi_\gamma(\beta)$ \ properly ends a \
$\Sigma_1(\nouse_\gamma)$--gap.
\ee
\end{corollary}

\begin{theorem}\label{thining} Let \ $\nouse$ \ be an iterable premouse and let \ $\gamma$ \ be an ordinal.
Then
\ $\boldface{\Sigma}{n}(\nouse_\gamma)\subseteq\boldface{\Sigma}{n}(\nouse)$ \ 
as pointclasses, for all \ $n\ge1$. 
\end{theorem}
\begin{proof} Corollary 2.20 of \cite{Crcm} directly implies this theorem.
\end{proof}

\begin{theorem}\label{preserve} If \ $\mouse$ \ is a real $1$--mouse, then for any premouse iterate \
$\mouse_\gamma$ \ we have that \ $\boldface{\Sigma}{n}(\mouse)=\boldface{\Sigma}{n}(\mouse_\gamma)$ \ 
as pointclasses, for all \ $n\ge1$. 
\end{theorem}
\begin{proof} Let \ $\mouse$ \  be a 1--mouse. We recall the definition of \ $\core=\core(\mouse)$, \ the core of \ $\mouse$.
\  Let \  $\mathcal{H} = \Hull_1^{\mouse}({\R\cup \omega\rho_{\mouse}} \cup \{p_{\mouse}\})$, \ and let   \ $\core$ \  be the
transitive collapse of \ $\mathcal{H}$. \ It follows that \ $\core$ \ is a real premouse and that \ $\mouse$ \ is a premouse
iterate \ of \ $\core$. \ Let \ $n\ge 1$. 
\begin{claim} $\boldface{\Sigma}{n}(\core)=\boldface{\Sigma}{n}(\mouse)$ \ as pointclasses.
\end{claim}
\begin{proof}[Proof of Claim] Theorem~\ref{thining} implies that \
$\boldface{\Sigma}{n}(\mouse)\subseteq\boldface{\Sigma}{n}(\core)$. \ Because \ $\core$ \ and \ $\H$ \ are isomorphic
structures, \
$\boldface{\Sigma}{n}(\C)=\boldface{\Sigma}{n}(\H)$. \ Since \ $\H$ \ is
\ $\Sigma_1$ \ definable (in parameters \ $\omega\rho_{\mouse}$ \ and \ $p_{\mouse}$)  over \ $\mouse$, \ we have that \
$\boldface{\Sigma}{n}(\H)\subseteq\boldface{\Sigma}{n}(\mouse)$. \ Therefore, \
$\boldface{\Sigma}{n}(\C)=\boldface{\Sigma}{n}(\mouse)$.
\end{proof}

Let \ $\mouse_\gamma$ \ be a premouse iterate of \ $\mouse$. \ Since \ $\mouse_\gamma$ \ is also a $1$--mouse with \
$\core(\mouse_\gamma)=\core(\mouse)$ \  (see Lemma 2.37 of \cite{Crcm}), the above Claim implies that \
$\boldface{\Sigma}{n}(\mouse)=\boldface{\Sigma}{n}(\mouse_\gamma)$  \ as pointclasses.
\end{proof}

\begin{theorem}\label{topless} Let \ $\mouse$ \ be a real 1--mouse. Then there is an ordinal \ $\gamma$ \
such that the premouse iterate \ $\mouse_\gamma$ \ is a proper initial segment of a real 1--mouse \ $\nouse$.
\end{theorem}
\begin{proof}
This follows immediately from Theorem 2.43 of \cite{Crcm}. 
\end{proof}

When \ $\mouse$ \ is a real 1--mouse, the next corollary shows that the question of
whether or not the external pointclass \  $\boldface{\Sigma}{n}(\mouse)$ \ has the scale property can be resolved, via Table
\ref{scaletable}, for each \ $n\ge 1$.
\begin{corollary}\label{toplesscor} Let \ $\mouse$ \ be a real 1--mouse. Then there exists an iterable real premouse \
$\nouse$ \ and an ordinal \ $\eta<\widehat{\OR}^{\,\nouse}$ \ such that \
$\boldface{\Sigma}{n}(\mouse)=\boldface{\Sigma}{n}(\nouse^\eta)$ \  as pointclasses for all \ $n\ge1$.
\end{corollary}
\begin{proof}
Let \ $\mouse$ \ be real 1--mouse. The proof of Lemma 5.4 of \cite{Crcm} implies that there is an ordinal \
$\gamma$ \ and a real 1--mouse \ $\nouse$ \ such that the premouse iterate \ $\mouse_\gamma$ \ is a proper initial
segment of \ $\nouse$. \ Note that \ $\pi_\gamma\colon \mouse \mapsigma{1} \mouse_\gamma$ \ preserves \
$\Sigma_1$--gaps by Theorem~\ref{gap_preseve}.  Let \ $\delta=\widehat{\OR}^{\,\mouse}$ \ and define \
$\pi(\delta)=\sup\left\{\pi(\gamma) :
\gamma<\widehat{\OR}^{\,\mouse}\right\}$. \ Note that \ $\pi(\delta)<\widehat{\OR}^{\,\nouse}$. \ It
follows that \ $\mouse_\gamma=\nouse^{\pi(\delta)}$. \ Hence, Theorem~\ref{preserve} implies that
\ $\boldface{\Sigma}{n}(\mouse)=\boldface{\Sigma}{n}\!\left(\nouse^{\pi(\delta)}\right)$ \
as pointclasses, for all \ $n\ge 1$. \  
Now, since \ $\pi(\delta)<\widehat{\OR}^{\,\nouse}$, \ we see that \ $\eta=\pi(\delta)$ \ is as desired. 
\end{proof}

\begin{theorem}\label{preserve2} If \ $\mouse$ \ is a real
mouse, then for any mouse iterate \ $\mouse_\gamma$ \ we have that 
\ $\boldface{\Sigma}{n}(\mouse)=\boldface{\Sigma}{n}(\mouse_\gamma)$ \ 
as pointclasses, for all \ $n\ge1$. 
\end{theorem}
\begin{proof} Lemma 2.19 of \cite{Cfsrm} implies this theorem.
\end{proof}

Using an argument similar to the one establishing Corollary~\ref{toplesscor} above, Theorem~\ref{preserve2} allows us to prove
our next corollary. Thus, when \ $\mouse$ \ is a real mouse  one  can also use Table~\ref{scaletable} to determine whether or
not the external pointclass \  $\boldface{\Sigma}{n}(\mouse)$ \ has the scale property, for \ $n\ge 1$.

\begin{corollary}\label{toplesscor2} Let \ $\mouse$ \ be a real mouse. Then there exists an iterable real premouse \
$\nouse$ \ and an ordinal \ $\eta<\widehat{\OR}^{\,\nouse}$ \ such that \
$\boldface{\Sigma}{n}(\mouse)=\boldface{\Sigma}{n}(\nouse^\eta)$ \  as pointclasses for all \ $n\ge1$.
\end{corollary}

\section{Scales of Minimal Complexity in $K(\mathbb{R})$}\label{minimalKr}
Since \ $\Kr$ \ is the union of real 1--mice, the development in Section~\ref{minimal} induces  a natural Levy hierarchy for
the sets of reals and the scales in \ $\Kr$. \ Before we identify the scales in \ $\Kr$ \ of minimal complexity,  we first show that
in \ $\Kr$ \ there is a close connection between the construction of scales and new \
$\Sigma_1$ \ truths about the reals. Assume that \ $\Kr\models\AD$. \   Given a scale \
$\seq{\le_i\, : i\in\omega}$, \ we shall denote this scale by \
$\seq{\le_i}$. \ Now let \ $A$ \ be a set of reals in \
$\Kr$ \ and suppose that \ $A$ \ has a scale \ $\seq{\le_i}$ \ in \ $\Kr$. \ By Theorem
5.5  of \cite{Crcm}, \ there is a real 1--mouse \ $\nouse\in\Kr$ \ such that \ $A\in\nouse$ \ and \ $\seq{\le_i}\in\nouse$. \
We will show that there is an ordinal \ $\gamma$ \ that begins a \ $\Sigma_1(\nouse)$--gap  in which \ $A$ \ has a \
$\boldface{\Sigma}{1}(\nouse^{\gamma})$ \ scale. Let \ $\eta$ \ be the least ordinal such that \ $A\in\nouse^{\eta+1}$. \
Let \ $[\alpha,\beta]$ \ be the \
$\Sigma_1(\nouse)$--gap containing \ $\eta$. \ Thus, \ $\alpha$ \ properly begins this gap. If \ $A$ \ is \
$\boldface{\Sigma}{1}(\nouse^{\alpha})$, \ then Theorem~\ref{alphascales} asserts that \ $A$ \ has a  \
$\boldface{\Sigma}{1}(\nouse^{\alpha})$
\ scale. If \ $A$ \ is not \ $\boldface{\Sigma}{1}(\nouse^{\alpha})$, \ it follows  that
\ $\beta$ \ must properly end this gap. To see this, suppose that \ $\nouse^\beta=\nouse$ \ and thus,  \ $A\in\nouse^\beta$ \
and \ $\seq{\le_i}\in\nouse^\beta$. \ Hence, by Wadge's Lemma every \ $\Pi_1(\nouse^\alpha)$ \ subset of
\ $\R\times\R$ \ has a uniformization in \ $\boldface{\Sigma}{1}(\nouse^{\beta})$. \ This contradicts Theorem
\ref{nouniform1}. Consequently, \ $\beta$ \ properly ends this gap and so, \ $\beta+1$ \ begins a new \ $\Sigma_1(\nouse)$--gap.
 Theorem~\ref{alphascales} implies that \ $A$ \ has a scale in \ $\boldface{\Sigma}{1}(\nouse^{\beta+1})$. \ Therefore, as
mentioned at the beginning of this paper, the construction of scales in \ $\Kr$ \ is closely tied to the verification of new \
$\Sigma_1$ \ truths in \ $\Kr$ \ about the reals. As a consequence of the above argument we have the following two theorems.

\begin{theorem} Suppose that \ $\nouse$ \ is an iterable real premouse such that \ $\nouse\models\AD$. \ If a set of
reals \ $A$ \ admits a scale in \ $\nouse$, \ then \ $A$ \ is \ $\Sigma_1(\nouse,\R)$.
\end{theorem}
\begin{theorem} Assume \ $\Kr\models\AD$. \ Any set of reals \ $A$ \ admits a
scale in \ $\Kr$ \ if and only if there is a real 1--mouse \ $\nouse$ \ such that \ $A\in
\boldface{\Sigma}{1}(\nouse^{\alpha})$ \ where
\ $\alpha$ \ begins a \ $\Sigma_1(\nouse)$--gap.
\end{theorem}

Furthermore, our method of defining scales in \ $\Kr$ \ produces 
scales of minimal complexity, as is established by Theorem~\ref{minnythm} below.   We shall now give a
precise definition of the ``complexity'' of a scale in \ $\Kr$. 
\begin{definition} Suppose that \ $\seq{\le_i}$ \ is a scale in \ $\Kr$. \ Let \ $\nouse$ \ be  a real
1--mouse such that \ $\seq{\le_i}\in\nouse$. \ Let \ $\seq{\gamma, k}$ \ be the lexicographically least so that
\ $\seq{\le_i}$ \ is \ $\boldface{\Sigma}{k}(\nouse^\gamma)$ \ or \ $\boldface{\Pi}{k}(\nouse^\gamma)$, \ and \
$\gamma<\widehat{\OR}^{\,\nouse}$. \ We shall say that
\ $\seq{\le_i}$ \ has \ $\nouse$--complexity \ $\seq{\gamma, k}$. \ Now let  \ $\seq{\le_i^*}$ \ be another scale in \
$\nouse$ \ and let \ $\seq{\gamma^*, k^*}$ \ be the 
$\nouse$--complexity of \ $\seq{\le_i^*}$. \ Then we shall say the $\nouse$--complexity of \ $\seq{\le_i}$ \ is less than or
equal to the $\nouse$--complexity of \ $\seq{\le_i^*}$, \ denoted by \ $\seq{\le_i}\preccurlyeq_\nouse\seq{\le_i^*}$, \ when \
$\seq{\gamma, k}$ \ is lexicographically less than or equal to \ $\seq{\gamma^*, k^*}$.
\end{definition}
Let \ $\nouse$ \ be any real 1-mouse containing the scale \ $\seq{\le_i}$ \ as an element. Lemma~{\minicomp} of \cite{Part1}
implies that the  $\nouse$--complexity of \ $\seq{\le_i}$ \ is always defined. 

\begin{definition} Let \ $\seq{\le_i}$ \ and \ $\seq{\le_i^*}$ \ be scales in \ $\Kr$. \ We shall say the
complexity of \ $\seq{\le_i}$ \ is less than or equal to the complexity of \ $\seq{\le_i^*}$, \ denoted by
\ $\seq{\le_i}\preccurlyeq\seq{\le_i^*}$, \ if and only if \ $\seq{\le_i}\preccurlyeq_\nouse\seq{\le_i^*}$ \ for some  1--mouse
\ $\nouse$.
\end{definition}
We note that the definition of the relation \ $\preccurlyeq$ \ is independent of the
1--mouse \ $\nouse$ \ (see Section~\ref{extensions}) and is a prewellordering on the scales in \ $\Kr$. \ 
\begin{theorem}\label{minnythm} Assume that \ $\Kr\models\AD$.
\ Let \ $\seq{\le_i^*}$ \  be a scale in \ $\Kr$ \ on a set of reals \ $A$. \ Then there is a scale \
$\seq{\le_i}$ \ on \ $A$ \ in \ $\Kr$ \ constructed as in Section~\ref{minimal} such that \
$\seq{\le_i}\preccurlyeq\seq{\le_i^*}$.
\end{theorem}
\begin{proof}
Let \ $\seq{\le_i^*}$ \ in \ $\Kr$ \ be a scale on a set of reals \ $A$. \ Let \ $\nouse$ \ be  a real 1--mouse such that \
$\seq{\le_i^*}\in\nouse$. \  Let \ $\eta$ \ be the least ordinal such that \
$\seq{\le_i^*}\in\nouse^{\eta+1}$. \ It follows that \ $A\in\nouse^{\eta+1}$. \  Let \ $[\alpha,\beta]$ \ be the \
$\Sigma_1(\nouse)$--gap containing \ $\eta$. \ Thus, \ $\alpha$ \ properly begins this gap. If \ $A$ \ is in \
$\boldface{\Sigma}{1}(\nouse^{\alpha})$, \ then Theorem~\ref{alphascales} asserts that \ $A$ \ has a  \
$\boldface{\Sigma}{1}(\nouse^{\alpha})$ \ scale \ $\seq{\le_i}$. \ Since \ $\seq{\le_i^*}\notin\nouse^\alpha$ \ it follows
that  \ $\seq{\le_i}\preccurlyeq_\nouse\seq{\le_i^*}$. \ Thus, if \ $A$ \ is \
$\boldface{\Sigma}{1}(\nouse^{\alpha})$, \ then the conclusion of the theorem follows. For the remainder of the
proof we shall assume that \ $A$ \ is not in \ $\boldface{\Sigma}{1}(\nouse^{\alpha})$. \ We focus on the
three cases: (1) $\eta=\alpha$, (2) $\alpha<\eta<\beta$ and (3) $\alpha<\eta=\beta$.

\smallskip
\noindent{\sc Case 1:} $\alpha=\eta$. \ Thus, \ $\seq{\le_i}\in\nouse^{\alpha+1}$ \ and \ $A\in\nouse^{\alpha+1}$. \ By Lemma~{\minicomp} of \cite{Part1}, there is a smallest natural number \ $n\ge 1$ \  such that \
$A\in\boldface{\Sigma}{n}(\nouse^{\alpha})\cup\boldface{\Pi}{n}(\nouse^{\alpha})$. 

\smallskip
{\sc Subcase 1.1:} {\sl $\alpha$  is uncollectible.\/}  Suppose that \ $n$ \ is odd. If  \ $A$ \ is in \
$\boldface{\Sigma}{n}(\nouse^{\alpha})$, \ then Theorem~\ref{pntclsequal} implies that \ $A$ \ has a
scale \ $\seq{\le_i}$ \ in \ $\boldface{\Sigma}{n}(\nouse^{\alpha})$. \ It follows that \
$\seq{\le_i}\preccurlyeq_\nouse\seq{\le_i^*}$. \ Suppose now that \ $A$ \ is not in \ $\boldface{\Sigma}{n}(\nouse^{\alpha})$
\ and thus, \ $A$ \ is in \ $\boldface{\Pi}{n}(\nouse^{\alpha})$. \ It follows that \ $\seq{\le_i^*}$ \ is not in \
$\boldface{\Pi}{n}(\nouse^{\alpha})$. \  Otherwise, since \ $A$
\ is not in \ $\boldface{\Sigma}{n}(\nouse^{\alpha})$, \ Wadge's Lemma would imply that \
$\boldface{\Pi}{n}(\nouse^{\alpha})$ \ has the scale property. In this case, however, \
$\boldface{\Sigma}{n}(\nouse^{\alpha})$ \ has the scale property and thus, \
$\seq{\le_i^*}$ \ cannot be in \ $\boldface{\Pi}{n}(\nouse^{\alpha})$ \ (see 4B.13 of \cite{Mosch}). We conclude that \
$A$ \ has a scale \ $\seq{\le_i}$ \ in \ $\boldface{\Sigma}{n+1}(\nouse^{\alpha})$ \ and  \
$\seq{\le_i}\preccurlyeq_\nouse\seq{\le_i^*}$. \ Similar
reasoning applies when \ $n$ \ is even. 

\smallskip
{\sc Subcase 1.2:} {\sl $\alpha$  is collectible.\/} This subcase is not possible, for suppose that \ $\alpha$ \ is
collectible. Since \ $\seq{\le_i^*}$ \ is in \ $\nouse^{\alpha+1}$ \ and \  $A$ \ is not in \
$\boldface{\Sigma}{1}(\nouse^{\alpha})$, \ Wadge's Lemma implies that the every \ $\Pi_1(\nouse^{\alpha})$ \ relation has a
uniformization in \ $\nouse^{\alpha+1}$, \ contradicting Theorem~\ref{nomore}.

\smallskip
\noindent{\sc Case 2:} $\alpha<\eta<\beta$. \ This case is not possible, for suppose that \ $\alpha<\eta<\beta$.  \ Since \
$\seq{\le_i^*}$ \ is in \ $\nouse^{\eta+1}$ \ and \  $A$ \ is not in \ $\boldface{\Sigma}{1}(\nouse^{\alpha})$, \ Wadge's Lemma
implies that the every \ $\Pi_1(\nouse^{\alpha})$ \ relation has a uniformization in \ $\boldface{\Sigma}{1}(\nouse^{\beta})$, \
contradicting Theorem~\ref{nouniform1}.

\smallskip
\noindent{\sc Case 3:} $\alpha<\eta=\beta$. \ Thus,  \
$\seq{\le_i}\in\nouse^{\beta+1}$, \ $A\in\nouse^{\beta+1}$, \ $A\notin\nouse^{\beta}$ \ and \ $\beta$ \ properly ends this
gap. Theorem~\ref{isamouse} implies that \ $\nouse^\beta$ \ is a real mouse and Theorem~\ref{newset} implies that \
$\m=m(\nouse^\beta)$ \ is defined. By Lemma~{\minicomp} of \cite{Part1}, there is a smallest natural number \ $n\ge 0$ \  such that \
$A\in\boldface{\Sigma}{\m+n}(\nouse^{\beta})\cup\boldface{\Pi}{\m+n}(\nouse^{\beta})$. \ 

\smallskip
{\sc Subcase 3.1:} {\sl $\beta$  is weak.\/}  Suppose that \ $n$ \ is even. If  \ $A$ \ is in \
$\boldface{\Sigma}{\m+n}(\nouse^{\beta})$, \ then Theorem~\ref{yesscale2} implies that \ $A$ \ has a
scale \ $\seq{\le_i}$ \ in \ $\boldface{\Sigma}{\m+n}(\nouse^{\beta})$. \ It follows that \
$\seq{\le_i}\preccurlyeq_\nouse\seq{\le_i^*}$. \ Suppose that \ $A$ \ is not in \ $\boldface{\Sigma}{\m+n}(\nouse^{\beta})$
\ and thus, \ $A$ \ is in \ $\boldface{\Pi}{\m+n}(\nouse^{\beta})$. \ Thus, \ $\seq{\le_i^*}$ \ is not in \
$\boldface{\Pi}{\m+n}(\nouse^{\beta})$. \  Otherwise, since \ $A$
\ is not in \ $\boldface{\Sigma}{\m+n}(\nouse^{\beta})$, \ Wadge's Lemma would imply that \
$\boldface{\Pi}{\m+n}(\nouse^{\beta})$ \ has the scale property. However, in this case, \
$\boldface{\Sigma}{\m+n}(\nouse^{\beta})$ \ has the scale property and hence, \
$\seq{\le_i^*}$ \ cannot be in \ $\boldface{\Pi}{\m+n}(\nouse^{\beta})$. \ Consequently, \ $A$ \ has a scale \  $\seq{\le_i}$ \ in \ $\boldface{\Sigma}{\m+n+1}(\nouse^{\beta})$ \ and  \ $\seq{\le_i}\preccurlyeq_\nouse\seq{\le_i^*}$. \ Similar reasoning applies when \ $n$ \ is odd. 

\smallskip
{\sc Subcase 3.2:} {\sl $\beta$  is strong.\/} This subcase is not possible, for suppose that \ $\beta$ \ is
strong. Since \ $\seq{\le_i^*}\in\nouse^{\beta+1}$, \ $A\in\nouse^{\beta+1}$ \  and \  $A$ \ is
not in \
$\boldface{\Sigma}{1}(\nouse^{\alpha})$, \ Wadge's Lemma implies that every \ $\Pi_1(\nouse^{\alpha})$ \
relation has a uniformization in \ $\nouse^{\beta+1}$, \ contradicting Theorem~\ref{nouniform2}.

This completes the proof.
\end{proof}

\section{Pointclass preserving premice}\label{questionQ} We now direct our attention to the question asked at the beginning
of this paper, namely:
\begin{question}[Q]Given an iterable real premouse  $\mouse$  and  $n\ge 1$,  when does the pointclass
$\boldface{\Sigma}{n}(\mouse)$  have the scale property?
\end{question}

Clearly, if an iterable real premouse \ $\mouse$ \ is a proper initial segment of an iterable real premouse, then the above
question can be addressed by referring to Table~\ref{scaletable}. Suppose now that \ $\mouse$ \ is not a proper initial
segment of another iterable real premouse. We know by Corollary 2.14(2) of \cite{Crcm} that the premouse iteration of \
$\mouse$ \ preserves the boldface pointclass \ $\boldface{\Sigma}{1}(\mouse)$, \ that is,
\ $\boldface{\Sigma}{1}(\mouse_\gamma)=\boldface{\Sigma}{1}(\mouse)$ \ for all ordinals \ $\gamma$. 
\ However, if \ $\mouse$ \ is also a real 1--mouse, then Theorem~\ref{preserve} asserts that the premouse
iteration of \ $\mouse$ \ preserves all of the boldface pointclasses, that is, \
$\boldface{\Sigma}{n}(\mouse_\gamma)=\boldface{\Sigma}{n}(\mouse)$ \ for all ordinals \ $\gamma$ \ and for all \ $n\ge 1$. 
Corollary~\ref{toplesscor} then implies that Question (Q) can be answered. Furthermore, when \ $\mouse$ \ is a real mouse, 
the fine structure of \ $\mouse$ \ can be used to prove that there is a mouse iterate \ $\mouse_\theta$
\ which is a proper initial segment of an iterable real premouse. Theorem~\ref{preserve2} asserts that mouse iteration
preserves the boldface pointclasses. Corollary~\ref{toplesscor2} thus implies that the question as to whether or not the
external pointclass \ $\boldface{\Sigma}{n}(\mouse)$ \ has the scale property can again be addressed. 

Let \ $\mouse$ \ be a  real mouse which is not a proper initial segment of an iterable real premouse. The above arguments show
that one can resolve Question (Q) by utilizing two fundamental attributes of \ $\mouse$:
\be
\item[(1)] $\mouse$ \ possesses a
specific fine structural property, and 
\item[(2)] $\mouse$ \ preserves the boldface pointclasses under mouse iteration.
\ee
This success inspires a general question. Suppose that \ $\mouse$ \ is merely an iterable real premouse that preserves the boldface
pointclasses under premouse iteration. {\sl Can it then be determined which, if any, of its external pointclasses have the scale
property?\/} It may be somewhat surprising to hear that the answer to this question is ``yes.'' If \ $\mouse$ \ is ``pointclass
preserving'',  we shall show that one can settle Question (Q) without presuming any specific fine structural conditions on \
$\mouse$.

\begin{definition}\label{presevedef} Let \ $\mouse$ \ be an iterable real premouse. We say that \ $\mouse$ \ is
{\it pointclass preserving\/} if, for every premouse iterate \ $\mouse_\gamma$ \ of \ $\mouse$, \ we have
that \ $\boldface{\Sigma}{n}(\mouse)=\boldface{\Sigma}{n}(\mouse_\gamma)$ \ as pointclasses, for each \ $n\ge 1$.
\end{definition}
\begin{rmk} When \ $\mouse$ \ is pointclass preserving, then every \
$\boldface{\Sigma}{n}(\mouse)$ \ set of reals \ $A$ \ is also \ $\boldface{\Sigma}{n}(\mouse_\gamma)$. \
Definition~\ref{presevedef} does not assert that the \ $\Sigma_n$ \ formula which defines \ $A$ \ over \ $\mouse$
\ is the same formula which defines \ $A$ \ over \ $\mouse_\gamma$. \ Definition~\ref{presevedef} only asserts
that there is {\it some\/} \ $\Sigma_n$ \ definition of \ $A$ \ over \ $\mouse_\gamma$ \ with parameters from \ $M_\gamma$.
\end{rmk}
For an iterable real premouse \ $\mouse$, \ Theorem~\ref{thining} asserts that \
$\boldface{\Sigma}{n}(\mouse_\gamma)\subseteq \boldface{\Sigma}{n}(\mouse)$ \ (as pointclasses) for any ordinal \
$\gamma$ \ and integer \ $n\ge 1$. \ Hence, premouse iteration does not produce any new definable sets of reals. Consequently, \
$\mouse$
\ is pointclass preserving if and only if for each ordinal \ $\gamma$ \ and \ $n\ge 1$ \ we have that \
$\boldface{\Sigma}{n}(\mouse)\subseteq\boldface{\Sigma}{n}(\mouse_\gamma)$. \ If \
$\mouse$ \ is pointclass preserving, then any premouse iterate of \ $\mouse$ \  does not ``lose'' any  \
$\boldface{\Sigma}{n}(\mouse)$ \ set of reals. For example, as noted above, real 1--mice are pointclass
preserving.

\begin{definition} Suppose that \ $\mouse$ \ and \ $\nouse$ \ are iterable
real premice.  Then 
\begin{align*}
\mouse \approx \nouse &\textup{ \ iff \ there exists a \ $\theta$ \ such that \ }
\mouse_\theta = \nouse_\theta,\\
\mouse \lesssim \nouse &\textup{ \ iff \ there exists a \ $\theta$ \ such that \ }
\mouse_\theta\textup{ \ is an initial segment of \ } \nouse_\theta,\\
\mouse < \nouse &\textup{ \ iff \ there exists a \ $\theta$ \ such that \ }
\mouse_\theta\textup{ \ is a proper initial segment of \ } \nouse_\theta.
\end{align*}
\end{definition} 
Recall that Definition~\ref{def:Rcomplete} identifies the notion of an $\R$--complete measure.
\begin{definition}\label{def:Rnormal} If \ $\nu$ is an $\R$--complete measure on \ $\kappa$ \ in \
$L[\nu](\R)$, \ then \ $L[\nu](\R)$ \ is said to be a  $\rho(\R)$--model with critical point \ $\kappa$.
\end{definition}
One can form repeated ultrapowers of a $\rho(\R)$--model \ $L[\nu](\R)$. \ If each such ultrapower is well-founded,
then we say that \ $L[\nu](\R)$ \ is {\it iterable\/}. 
The next theorem shows that if an iterable real premouse \ $\mouse$ \
is ``larger'' than all real 1--mice, then there is an iterable
$\rho(\R)$--model which contains an iterate of \ $\mouse$ \ as a proper initial segment. If \ $\mouse$ \ is pointclass
preserving,  then this theorem will allow us to determine if any of the pointclasses \ $\boldface{\Sigma}{n}(\mouse)$ \ have the scale
property.
\begin{theorem}\label{highestmouse} Suppose that \ $\mouse$ \
is an iterable real premouse such that \ $\nouse\lesssim\mouse$ \ for all real 1--mice \ $\nouse$. \ Then there
exists an iterable $\rho(\R)$--model \ $L[\nu](\R)$ \ with critical point \ $\lambda>\kappa$ \ and an ordinal \
$\theta$ \ such that the premouse iterate \ $\mouse_\theta$ \ is an initial segment of \
$L[\nu](\R)$. \ In addition, \ $\pow(\R)\cap L[\nu](\R)\subseteq\Kr$.
\end{theorem}

\begin{proof} Let \ $\mouse=(M,\R,\kappa,\mu)$. \ It follows from Theorem~\ref{topless} that  \ $\nouse<\mouse$ \
for all real 1--mice \ $\nouse$. \ The proof of Lemma 5.4 of \cite{Crcm} \ can be used to show that every subset
of \ $\kappa$ \ in \ $\Kr$ \ is also in some iterate of \ $\mouse$. \ Hence, Lemma 2.11 of \cite{Crcm} implies that every
such subset of \ $\kappa$ \ is in \ $\mouse$. \ A similar argument shows that any \ $\kappa$--sequence 
of subsets of \ $\kappa$ \ in \ $\Kr$ \ is also in \ $\mouse$. \ Therefore, \ $\mu$ \ is a \ $\Kr$--measure on \
$\kappa>\Theta^{\Kr}$ (see Definition 2.1 of \cite{Cnot}). Without loss of generality, one can assume that the
measure \ $\mu$ \ is countably complete.\footnote{If not, then apply the argument to \
$\mouse_\lambda$ \ where \ $\lambda$ \ is a sufficiently large regular cardinal.} By \ $\DC$ \ and countable
completeness, it follows that the ultrapower \ ${^\kappa}\Kr/\mu$ \ is well-founded. The proof of Corollary 2.14 of \cite{Cnot}
implies that \ $(\Kr,\mu)$ \ is really good on \ $\kappa$ \ (see \cite[Definition 2.9]{Cnot}).  Again, by \ $\DC$ \ and countable
completeness, \ $(\Kr,\mu)$ \ is weakly iterable (see \cite[Section 2]{Cnot}). Lemma 4.3 of \cite{Cnot} now implies the existence of
the desired
$\rho(\R)$--model \ $L[\nu](\R)$.
\end{proof}

\begin{remark}
Let \ $\mouse$ \ and \ $L[\nu](\R)$ \ be as in the statement of Theorem~\ref{highestmouse}. If \ $\Kr\models\AD$, \
then it follows that \ $L[\nu](\R)\models\AD$ \ and Lemma 2.11 of \cite{Crcm} implies that \
 $\mouse\models\AD$.
\end{remark}

\begin{corollary} Let \ $\mouse$ \ be a pointclass preserving premouse such
that \ $\nouse\lesssim\mouse$ \ for all real 1--mice \ $\nouse$. \ For all \ $n\ge 1$ \ and all sets of reals \ $A$, \
if \  $A$ \ is \ $\boldface{\Sigma}{n}(\mouse)$ \ then \ $A\in\mouse$.
\end{corollary}
\begin{proof} Let \ $n\ge 1$ \ and suppose that \ $A$ \ is a set of reals in \ $\boldface{\Sigma}{n}(\mouse)$. \ Let \ $\theta$ \ be as in
Theorem~\ref{highestmouse}. \ Since \ $\mouse$ \ is pointclass preserving, we have that  \
$\boldface{\Sigma}{n}(\mouse_\theta)=\boldface{\Sigma}{n}(\mouse)$ \ as pointclasses.  Theorem
\ref{highestmouse} implies that \ $A\in\Kr$. \ Since \ $\nouse\lesssim\mouse$ \ for all real 1--mice \ $\nouse$, \ it follows that \
$\nouse<\mouse$ \ for all real 1--mice \ $\nouse$. \ Thus, \ $A$
\ must be an element of \ $\mouse$.
\end{proof}

\begin{corollary} Assume that \ $\Kr\models\AD$. \  Let \ $\mouse$ \ be a pointclass preserving
premouse. Suppose that \ $\nouse\lesssim\mouse$ \ for all real 1--mice \ $\nouse$. \ Then, 
\be
\item $\Sigma_1(\mouse,\R)$ \ has the scale property, and
\item $\boldface{\Sigma}{n}(\mouse)$ \ and \ $\boldface{\Pi}{n}(\mouse)$ \ do not have the scale property
for any \ $n\ge 1$.
\ee
\end{corollary}
\begin{proof} Theorem~\ref{highestmouse} implies that \ $\mouse\models\AD$ \ and thus, \ $\Sigma_1(\mouse,\R)$ \
has the scale property by Corollary~\ref{relscales}.   Theorem 5.17 of \cite{Crcm} implies that there is a set of reals \ $A$ \ in \ 
$\Kr$ \ which has no scale in \ $\Kr$. \ Theorem~\ref{highestmouse} implies that \ $A\in\mouse$. \ Using \
$A$ \ as a constant, it follows that \ $\boldface{\Sigma}{n}(\mouse)$ \ does not have the scale property for any \
$n\ge 1$.
\end{proof}
\begin{theorem}Assume that \ $\Kr\models\AD$. \ Let \ $\mouse$ \ be a pointclass preserving premouse. Suppose
that \ $\mouse< \nouse$ \ for some real 1--mice \ $\nouse$. \ Then, for \ $n\ge 1$ \ one can determine whether
or not the pointclass \ $\boldface{\Sigma}{n}(\mouse)$, \ or \ $\boldface{\Pi}{n}(\mouse)$, \ has the scale
property.
\end{theorem}
\begin{proof} Suppose that \ $\mouse$ \ is a pointclass preserving premouse such that \ $\mouse< \nouse$ \ for some real
1--mice \ $\nouse$. \ Let \ $\theta$ \ be so large that \ $\mouse_\theta$ \ is an initial segment of \ $\nouse_\theta$. \ Thus,
\ $\boldface{\Sigma}{n}(\mouse)=\boldface{\Sigma}{n}(\mouse_\theta)$ \ as pointclasses.  Since $\mouse_\theta$ \ is a proper
initial segment of \ $\nouse_\theta$, \ let \ $\gamma<\widehat{\OR}^{\nouse_\theta}$ \ be such that \
$\mouse_\theta=\nouse_\theta^\gamma$. \ Let \ $[\alpha,\beta]$ \ be the \ $\Sigma_1(\nouse_\theta)$--gap containing \ $\gamma$.
\ Now Table~\ref{scaletable} can be used to determine whether or not \ $\boldface{\Sigma}{n}(\mouse)$ \ has the scale property.
\end{proof}


\providecommand{\bysame}{\leavevmode\hbox to3em{\hrulefill}\thinspace}
\providecommand{\href}[2]{#2}

\end{document}